\documentclass[preprint,12pt]{elsarticle}

\usepackage{amsbsy,amssymb,amsmath,graphicx}
\usepackage[rightcaption]{sidecap}
\usepackage{multirow}
\usepackage{graphics}
\usepackage{color}

\begin{document}

\begin{frontmatter}



\title{On a local Fourier analysis for overlapping block smoothers on triangular grids}


\author{C. Rodrigo\corref{cor1}}
\author{F.~J. Gaspar}
\author{F.J. Lisbona}
\cortext[cor1]{Corresponding author. Tel:~+34976762148. \emph{E-mail~address}:~carmenr@unizar.es}
\address{Applied Mathematics Department, University of Zaragoza, Spain}

\begin{abstract}
A general local Fourier analysis for overlapping block smoothers on triangular grids is presented. This analysis is explained in a general form for its application to problems with different discretizations. This tool is demonstrated for two different problems: a stabilized linear finite element discretization of Stokes equations and an edge-based discretization of the curl-curl operator by lowest-order N\'{e}d\'{e}lec finite element method.
In this latter, special Fourier modes have to be considered in order to perform the analysis.
Numerical results comparing two- and three-grid convergence factors predicted by the local Fourier analysis to real asymptotic convergence factors are presented to confirm the predictions of the analysis and show their usefulness.
\end{abstract}

\begin{keyword}
Multigrid \sep triangular grids \sep overlapping block smoothers \sep box-relaxation \sep Vanka smoothers
\sep local Fourier analysis \sep saddle point type problems \sep Stokes
\sep N\'{e}d\'{e}lec


\end{keyword}

\end{frontmatter}

\section{Introduction}\label{sec:Intro}

As is well-known, multigrid methods \cite{Bra77, Hackb} are among the most powerful techniques for the efficient resolution of the large systems of equations arising from the discretization of partial differential equations. Since the 70's, when these methods were developed, they have become very popular among the scientific community. They have the nice property of requiring a computational work of the order of the number of unknowns of the problem, at least for elliptic problems. Besides, they have also been applied to more complicated problems, for example see \cite{TOS01}, providing very good results.

The efficiency and the robustness of a multigrid method is essentially influenced by the smoothing algorithm. We want to study the class of multiplicative Schwarz smoothers. Basically, they can be described as an overlapping block Gauss-Seidel method, where a number
of small linear systems of equations has to be solved in each smoothing step. This type of smoothers is characterized by its ability to deal with
saddle point problems and equations where the terms grad-div or curl-curl dominate. A particular case of such relaxation is the so-called Vanka smoother, introduced in~\cite{vanka} for solving the staggered finite difference discretization of the Navier-Stokes equations.

Local Fourier analysis (LFA, or local mode analysis) is a commonly used approach for analyzing the convergence properties of geometric multigrid methods.
In this analysis an infinite regular grid is considered and boundary conditions are not taken into account. LFA was introduced by Brandt in \cite{Bra77} and afterward extended in \cite{Bra94}. A good introduction can be found in the paper by St\"{u}ben and Trottenberg \cite{Stu_Tro} and in the books by Wesseling \cite{Wess}, Trottenberg et al. \cite{TOS01}, and Wienands and Joppich \cite{Wie01}. LFA was generalized to triangular grids in \cite{Gaspar_LFA}, for discretizations based on linear finite element methods. Afterwards, this generalization has been extended to systems of partial differential equations
\cite{mma2010, elasticity} and to high-order finite element discretizations \cite{cuadraticos}.

To our knowledge, there are only few papers dealing with local Fourier analysis for overlapping smoothers, all of them for discretizations on rectangular grids.  This analysis was performed in \cite{Siv-91} for the staggered finite-difference discretizations of the Stokes equations, and in \cite{molenaar-91} for a mixed finite element discretization of the Laplace equation. In~\cite{vandewalle2008}, an LFA to analyze an additive Schwarz smoother for a curl-curl model problem is proposed. A multicolored version was considered in the way that the corresponding analysis does not consider the special techniques to study the standard overlapping smoothers. In \cite{ThesisRodrigo} an LFA for overlapping block smoothers on triangular grids is presented. This tool was applied to linear finite element discretizations for poroelasticity problems. Later, in  \cite{Oost-2010}, the analysis for such overlapping block smoothers is performed on rectangular grids for finite element discretizations of the grad-div, curl-curl and Stokes equations. Here, we present and extend this analysis to general discretizations on triangular grids, including some special techniques for the case of edge-based discretizations. Two model problems are chosen to show this analysis, but we keep in mind that it can be carried over to a variety of other problems and other overlapping smoothers. The considered problems
are the discretization by stabilized linear finite elements of the Stokes problem, and the low-order N\'{e}d\'{e}lec's edge elements for the curl-curl equation. Regarding the Stokes system, we perform an exhaustive two-grid local Fourier analysis for the full- and diagonal (much cheaper) versions of the overlapping block smoother. Apart from this, also a three-grid analysis is developed in order to obtain more insight. This analysis has to deal with the difficulties inherent to the treatment of this kind of smoothers and also it has to take into account the corresponding extension of LFA to triangular grids. Apart from these difficulties, for the second model problem we need to extend this analysis to edge-based discretizations on which different stencils appear depending on the type of edge. This makes necessary the introduction of special Fourier modes in order to perform the analysis. Also for this test we perform a two- and three-grid LFA to analyze the differences in the performance of W- and V-cycles.

The structure of the paper is as follows.
In Section~\ref{sec:vanka}, a general description of the class of overlapping block smoothers is done, together with the development of a suitable local Fourier analysis for this type of relaxation procedures. Two- and three-grid local Fourier analysis are performed. Sections~\ref{sec:saddle} and~\ref{sec:vector} are focused on a detailed description of the local Fourier analysis of two particular overlapping block smoothers for the solution of two different model problems. More concretely, in Section~\ref{sec:saddle} a stabilized linear finite element discretization of the Stokes equations is considered and in Section~\ref{sec:vector} we deal with an edge-based discretization of a curl-curl problem by using low order N\'{e}d\'{e}lec finite elements.
Finally, in Section~\ref{sec:conclusions} some conclusions are drawn.

\section{Local Fourier analysis for overlapping block smoothers}\label{sec:vanka}

\subsection{Description of the smoother}\label{sec:smoother}

Point-wise iterative methods can be generalized to block-wise iterative schemes by updating a set of unknowns at each time, instead of only one. To this end, the grid is split into blocks and the equations corresponding to the grid-points in each block are simultaneously solved as a system of equations. Block-wise schemes become very attractive when anisotropies appear, especially when they are combined with a problem-dependent ordering of the blocks, since point-wise relaxation techniques lose their smoothing property.
Many arbitrary splittings of the mesh can be considered. For example, it is possible to allow the blocks to overlap, what gives rise to the class of overlapping block iterations, where smaller local problems are solved and combined via a multiplicative Schwarz method. They were introduced by Vanka in \cite{vanka}, and in~\cite{Zulehner} a theoretical basis for this approach was provided.
These smoothers, also known as coupled or box-relaxation, consist of decomposing the mesh into small subdomains and treating them separately.  Therefore, one relaxation step consists of a loop over all subdomains, solving for each one the system arising from the corresponding equations. Next, we give a more detailed description of the iterative method. We consider a linear system of equations $A_h\,u_h = f_h,$ which, in our case, arises from the discretization of a PDE problem. Vector $u_h$ is composed of unknowns corresponding to $m$ different variables. More concretely, $N_i$ unknowns of each variable $i$ are considered. Let $B$ be the subset of unknowns involved in an arbitrary block, that is, $B = \{u^1_{k_1(1)},\ldots, u^1_{k_1(n_1)},\ldots,u^m_{k_m(1)},\ldots,u^m_{k_m(n_m)}\},$ where $k_i(1),\ldots,k_i(n_i)$ are the global indexes of the $n_i$ unknowns corresponding to variable $i$. In order to obtain the matrix $A_h^B$ of the system to solve associated with block $B,$ we introduce the matrix $V_B$ representing the projection operator from the vector of all unknowns to the vector of the unknowns involved in the block, as the following block-diagonal matrix
\begin{equation}\label{projection_matrix}
V_B = \left(\begin{array}{ccc} V_B^1 &  & \\  & \ddots &  \\  &  & V_B^m \end{array}\right).
\end{equation}
Here, each block $V_B^i$ is a $(n_i\times N_i)-$matrix, whose $j$th-row is the $k_i(j)$th-row of the identity matrix of order $N_i.$
In this way, matrix $A_h^B$ can be defined as
\begin{equation}\label{block_matrix}
A_h^B = V_B\,A_h\, V_B^T.
\end{equation}
Therefore, this type of smoother results in a multiplicative Schwarz method with iteration matrix
\begin{equation}
\prod_{B=1}^{NB}(I-V_B^T\,(A_h^B)^{-1}\,V_B\,A_h),
\end{equation}
where $NB$ is the number of blocks or small systems to be solved in a relaxation step of the iterative method. Very often in practice, instead of solving the local problems exactly, one can replace $A_h^B$ with an approximation ${\widetilde A}_h^B$, obtaining the so-called multiplicative Schwarz method with inexact local solver, with iteration matrix given by
\begin{equation}
\prod_{B=1}^{NB}(I-V_B^T\,(\widetilde{A}_h^B)^{-1}\,V_B\,A_h).
\end{equation}
Therefore, many variants of box-relaxation can be considered. They can differ in the choice of the subdomains which are solved simultaneously, and in the way in which the local systems to be solved are built. Notice that the different subdomains can also be visited in different orderings, for instance, they can also be treated with some patterning scheme, yielding to a multicolored version of these relaxation schemes. All this makes wider the variety of this type of smoothers.\\
Regarding the choice of the subdomains, the first one coming into mind is the so-called cell-based box smoother. As indicated by its name, the considered subdomains are the elements of the triangulation.
Then, the unknowns located at each cell are simultaneously relaxed.
Another variant of this kind of smoothers is the point-box smoother, which consists of performing a loop over all grid-points, considering some of the unknowns appearing in the stencil of each node of the grid, and solving simultaneously the system arising from the equations associated with these points.

While much literature can be found presenting numerical studies of different overlapping block smoothers, there are only few papers that deal with theoretical aspects of these relaxation schemes~\cite{Siv-91, molenaar-91, ThesisRodrigo, Oost-2010, Zulehner, manservisi}. Here, we want to describe in a general setting how to perform local Fourier analysis for these smoothers on the general case of triangular grids (notice that rectangular and hexagonal grids are included in this setting).

\subsection{Fourier analysis}\label{sec:lfa}

Overlapping block smoothers require a special strategy to carry out their local Fourier analysis, since a classical analysis fails for this class of smoothers. The distinction with respect to other smoothers is that this relaxation updates some variables more than once, due to the overlapping of the local subdomains which are simultaneously solved. This fact has to be taken into account in the analysis because it causes that, in addition to the initial and final errors, some intermediate errors appear. Next, this non-standard local Fourier analysis is described for general discretizations on triangular grids. With this aim, the ideas about the recently introduced LFA on triangular grids~\cite{Gaspar_LFA, Rodrigo} have to be taken into account. The key fact for this extension is to consider an expression of the Fourier transform in new coordinate systems in space and frequency variables. We establish a non-orthogonal unitary basis of ${\mathbb R}^2,$ $\{{\mathbf e}_1,{\mathbf e}_2\},$ fitting the geometry of the given mesh, and the basis corresponding to the frequency space, $\{{\mathbf e}_1',{\mathbf e}_2'\},$ is taken as its reciprocal basis; that is, the vectors of the bases satisfy $({\mathbf e}_i,{\mathbf e}'_j)=\delta_{ij}, \, 1\le i,j\le 2,$ see~\cite{Gaspar_LFA, Rodrigo} for more details. Notice that in the case of a cartesian grid, the reciprocal basis is the same as the original one and then it arises as a particular case in this general framework.

To describe the analysis for the overlapping block smothers we follow the methodology developed in~\cite{Siv-91}, although in this work we present this analysis in a very general framework, with an arbitrary number of variables. For simplicity in the presentation, we restrict ourselves to the case of blocks of the same size, $q,$ which is the most common approach in practice and also a natural setting in the LFA framework. This means that for each block, the $(q\times q)-$system composed of the equations corresponding to such unknowns has to be solved. Notice that different variables can be involved in the block which is going to be simultaneously updated, as for example in Stokes equations where some velocity and pressure variables are updated together as we will see in Section~\ref{sec:saddle}. Assuming that unknowns associated with $m$ variables are involved, the following system, in terms of corrections and residuals, arises for the block number~$B,$
\begin{equation}\label{vanka_system}
A_h^{B} \, \left(\begin{array}{c}{\boldsymbol \delta} {\mathbf u}^1 \\ \vdots \\ {\boldsymbol \delta} {\mathbf u}^m\end{array}\right) = \left(\begin{array}{c}{\mathbf r}^1 \\ \vdots \\ {\mathbf r}^m\end{array}\right),
\end{equation}
where ${\boldsymbol \delta} {\mathbf u}^i$ is the vector of the corrections corresponding to the unknowns associated with variable $i$ involved in such block. Given an arbitrary unknown of the block, corresponding to variable $i,$ located at node $(k,l),$ it is considered that $\delta u^i_{k,l} = u_h^{i,j+(n+1)/s}({\mathbf x}_{k,l})-u_h^{i,j+n/s}({\mathbf x}_{k,l}),$ and\\ $r^i_{k,l}~= f_h^i({\mathbf x}_{k,l}) - A_hu^{i,j+n/s}_h({\mathbf x}_{k,l}),$ where $u_h^{i,j+n/s}({\mathbf x}_{k,l})$ denotes the approximation of $u^i_h$ computed in the $j-$iteration and which has been already updated $n$ times in the current relaxation step. Here, it is considered that the corresponding unknown can be updated up to $s$ times, and ${\mathbf x}_{k,l}$ denotes the grid-point associated with numbering $(k,l).$
This system can be written in terms of the errors, but to this end, we must take into account that some intermediate errors appear. Thus, we denote as $e_h^{i,j}({\mathbf x})$ the initial error of variable $i$ at $j-$iteration, $e_h^{i,j+1}({\mathbf x})$ the initial error of variable $i$ at $(j+1)-$iteration, and $e_h^{i,j+n/s}({\mathbf x}),\, n=1,\ldots,s-1$ the error obtained for variable $i$ after the unknown in node ${\mathbf x}$ has been updated $n-$times in the current iteration.
As stated the local Fourier analysis assumptions, errors can be written as a formal linear combination of Fourier modes.
An important aspect to take into account in the analysis is whether the Fourier modes are eigenfunctions of the smoothing operator. This is fulfilled for most standard relaxation operators, but there are some exceptions as for example pattern relaxations. In the case of the overlapping block smoothers treated here, this is not straightforward. A rigorous and detailed proof of this fact can be found in the literature in the paper by Maclachlan et al~\cite{Oost-2010}, where the authors prove that Fourier modes are also eigenfunctions of any coupled relaxation on rectangular grids. This allows to perform the analysis of overlapping block smoothers using classical LFA techniques on rectangular grids. The extension of this idea to triangular grids does not represent any substantial change in the proof presented in~\cite{Oost-2010} for the cartesian grid case, so we refer to the reader to this proof for a detailed understanding of the issue.
Then, without loss of generality, let the initial and the fully corrected errors be given by a single Fourier mode, that is
\begin{eqnarray}\label{errors_fourier}
e_h^{i,j}({\mathbf x}) &=& \alpha^{(0)}_i(\boldsymbol \theta) e^{\imath\boldsymbol\theta {\mathbf x}},\nonumber\\ [-1ex]
\\
e_h^{i,j+1}({\mathbf x}) &=& \alpha^{(s)}_i(\boldsymbol \theta) e^{\imath\boldsymbol\theta {\mathbf x}},\nonumber
\end{eqnarray}
where the coefficient of the Fourier mode is characterized by a super-index indicating the number of times that the unknown has been already corrected and by a subindex denoting the corresponding variable. Following this criterion, we denote the intermediate errors as
$$e_h^{i,j+n/s}({\mathbf x}) = \alpha^{(n)}_i(\boldsymbol \theta) e^{\imath\boldsymbol\theta {\mathbf x}},\; n=1,\ldots,s-1.$$
The aim is therefore to find the relation between the initial and the fully corrected errors for the involved variables, which is given by $$\left(\begin{array}{c}\alpha^{(s)}_1(\boldsymbol \theta)\\ \vdots \\ \alpha^{(s)}_m(\boldsymbol \theta)\end{array}\right) = \widetilde{S}_h(\boldsymbol \theta) \left(\begin{array}{c}\alpha^{(0)}_1(\boldsymbol \theta)\\ \vdots \\ \alpha^{(0)}_m(\boldsymbol \theta)\end{array}\right).$$ To this end, system (\ref{vanka_system}) can be written in terms of the errors, and consequently in terms of the coefficients $\alpha^{(n)}_i(\boldsymbol \theta),$ by taking into account that the corrections can also be written as $\delta u_{k,l}^i = e_h^{i,j+(n+1)/s}({\mathbf x})-e_h^{i,j+n/s}({\mathbf x}).$
Next, the resulting system, can be rearranged into a system of equations for the updated Fourier coefficients: $\alpha^{(1)}_i(\boldsymbol \theta),\ldots, \alpha^{(s)}_i(\boldsymbol \theta),$ leaving in the right-hand side the terms in $\alpha^{(0)}_i(\boldsymbol \theta),$ corresponding to the non-updated components. Thus, we obtain a system of the form
\begin{equation}\label{systemPQ}
P\;\left( \begin{array}{c}\alpha^{(1)}_1(\boldsymbol \theta)\\ \vdots\\ \alpha^{(1)}_m(\boldsymbol \theta)\\ \vdots \\ \alpha^{(s)}_1(\boldsymbol \theta)\\ \vdots \\ \alpha^{(s)}_m(\boldsymbol \theta) \end{array} \right) = \;Q\;\left(\begin{array}{c} \alpha^{(0)}_1(\boldsymbol \theta)\\ \vdots \\ \alpha^{(0)}_m(\boldsymbol \theta)\end{array}\right),
\end{equation}
where $P$ is a $((s\,m)\times (s\,m))-$matrix and Q is $(s\,m)\times m.$
As system (\ref{systemPQ}) can be written as
\begin{equation}\label{systemPQ_2}
\left( \begin{array}{c}\alpha^{(1)}_1(\boldsymbol \theta)\\ \vdots\\ \alpha^{(1)}_m(\boldsymbol \theta)\\ \vdots \\ \alpha^{(s)}_1(\boldsymbol \theta)\\ \vdots \\ \alpha^{(s)}_m(\boldsymbol \theta) \end{array} \right) = \;(P^{-1}\,Q)\; \left(\begin{array}{c} \alpha^{(0)}_1(\boldsymbol \theta)\\ \vdots \\ \alpha^{(0)}_m(\boldsymbol \theta)\end{array}\right),
\end{equation}
the last $(m\times m)-$block-entry of the $((s\,m)\times m)-$matrix $P^{-1}\,Q$ results to be the amplification matrix for the complete sweep, since it relates the Fourier coefficients of the fully corrected and the initial errors, that is $$\left(\begin{array}{c} \alpha^{(s)}_1(\boldsymbol \theta)\\ \vdots \\ \alpha^{(s)}_m(\boldsymbol \theta)\end{array}\right) = (P^{-1}\,Q)_{\{(s-1)m+1:s\,m, 1:m\}} \; \left(\begin{array}{c} \alpha^{(0)}_1(\boldsymbol \theta)\\ \vdots \\ \alpha^{(0)}_m(\boldsymbol \theta)\end{array}\right),$$ what means that $\widetilde{S}_h(\boldsymbol \theta) = (P^{-1}\,Q)_{\{(s-1)m+1:s\,m, 1:m\}}.$
Notice that if scalar problems are considered ($m=1$), then the Fourier symbol $\widetilde{S}_h(\boldsymbol \theta)$ consists on a single number instead of a matrix.
Finally, by considering the obtained Fourier representation for the relaxation process, the smoothing as well as the k-grid Fourier analysis for this smoother can then be performed as usual. Next, a brief description of these well-known techniques is presented, but for a deeper insight we refer the reader to~\cite{TOS01, Wie01}. \\
One can perform a local Fourier smoothing analysis by considering the corresponding Fourier domain representation of the smoothing operator $\widetilde{S}_h({\boldsymbol \theta})$ on the high frequencies, that is, on the subset ${\boldsymbol \Theta}_h\backslash (-\pi/2,\pi/2]^2$, with ${\boldsymbol \Theta}_h = (-\pi,\pi]^2$, when standard coarsening is considered.
This gives rise to the computation of the so-called smoothing factor, which is defined as
\begin{equation}\label{smoothing_factor}
\mu = \sup_{{\boldsymbol \Theta}_h \backslash (-\pi/2,\pi/2]^2} \rho(\widetilde{S}_h({\boldsymbol \theta})).
\end{equation}
However, in order to investigate the interplay between relaxation and coarse--grid correction, which is crucial for an efficient multigrid method, it is necessary to perform at least a two--grid analysis which takes into account the effect of transfer operators. For this purpose, one needs to consider the error propagation operator from the two--grid method, that is,
$$M_h^{2h} = S_h^{\nu_2}(I_h-I_{2h}^hA_{2h}^{-1}I_h^{2h}A_h)S_h^{\nu_1},$$
where $S_h$ is the smoothing procedure and the coarse--grid correction operator is composed of the discrete operators on the fine and coarse grids, $A_h$ and $A_{2h}$, respectively, and the inter--grid transfer operators: restriction, $I_h^{2h}$ and prolongation $I_{2h}^h$.
The two–-grid analysis is the basis for the classical asymptotic multigrid convergence
estimates, and the spectral radius of the two--grid operator, $\rho(M^{2h}_h)$, indicates
the asymptotic convergence factor of the two-grid method.
To estimate this value, the crucial observation is that the coarse--grid correction operator, as well as the smoother, leave the so-called spaces of $2h$-harmonics, ${\mathcal F}^4({\boldsymbol \theta}^{00})$, invariant. These subspaces are given by
$${\mathcal F}^4({\boldsymbol \theta}^{00}) = span\{\phi_h({\boldsymbol \theta}^{\alpha_1 \alpha_2},{\mathbf x})| \alpha_1,\alpha_2\in\{0,1\}\}, \; \hbox{with} \; {\boldsymbol \theta}^{00}\in(-\pi/2,\pi/2]^2,$$
and where ${\boldsymbol \theta}^{\alpha_1 \alpha_2}={\boldsymbol \theta}^{00}-(\alpha_1 \hbox{sign}(\theta^{00}_1)\pi,\alpha_2 \hbox{sign}(\theta^{00}_2)\pi)$.
For this reason, $M_h^{2h}$ is equivalent to a block--diagonal matrix consisting of blocks denoted by $\widetilde{M}_h^{2h}({\boldsymbol \theta}^{00}) = M_h^{2h}|_{{\mathcal F}^4({\boldsymbol \theta}^{00})}$, that is, its Fourier domain representation. In this way, one can determine the spectral radius $\rho(M_h^{2h})$ by calculating the spectral radius of these smaller matrices, that is:
\begin{equation}\label{rho}
\rho_{2g} = \rho(M_h^{2h}) = \displaystyle\sup_{{\boldsymbol \theta}^{00}\in(-\pi/2,\pi/2]^2}\rho(\widetilde{M}_h^{2h}({\boldsymbol \theta}^{00})).
\end{equation}
Furthermore, a deeper insight into the performance of multigrid can be obtained by considering a $k-$grid analysis ($k>2$), see~\cite{Wie01}. For example, the influence of the type of cycle (V-cycle versus W-cycle) and/or the number of pre- and post-smoothing steps can not be accurately predicted by a classical Fourier two--grid analysis. Due to the recursive definition of the k-grid operator, the previously introduced two--grid analysis can be extended to a k-grid analysis, and in particular to a three-grid analysis. Analogously to the case of the two-grid cycle, the expression giving the transformation of the error by a three-grid cycle is $e_h^{m+1} = M_h^{4h} e_h^m,$ with $$M_h^{4h} = S_h^{\nu_2} C_h^{4h} S_h^{\nu_1}= S_h^{\nu_2} (I_h-I_{2h}^h(I_{2h}-(M_{2h}^{4h})^{\gamma})(A_{2h})^{-1}I_h^{2h}A_h)S_h^{\nu_1},$$
being $M_{2h}^{4h}$ the two-grid operator between the two coarse grids, defined as
$$M_{2h}^{4h} = S_{2h}^{\nu_2}(I_{2h}-I_{4h}^{2h}(A_{4h})^{-1}I_{2h}^{4h}A_{2h})S_{2h}^{\nu_1},$$
and $\gamma$ the number of two-grid iterations. As well as for the two--grid analysis was crucial the decomposition of the Fourier space into the subspaces of $2h-$harmonics, for the three-grid analysis it has to be taken into account that not only in the transition from the finest to the second grid but also in the transition from the second grid to the coarsest grid four frequencies are coupled. Taking this into account, we can define appropriate subspaces of $4h-$harmonics, generated for the Fourier components associated with the sixteen frequencies coupled on the coarsest grid
\begin{equation}
{\cal F}^{16}({\boldsymbol \theta}^{00}) = span\{\varphi_h({\boldsymbol \theta}^{ij}_{nm}),\; i,j,n,m \in\{0,1\}\},
\end{equation}
with ${\boldsymbol \theta}^{00} = (\theta_1^{00},\theta_2^{00}) \in{\boldsymbol \Theta}_{4h}=(-\pi/4,\pi/4]\times(-\pi/4,\pi/4],$ and
where
\begin{eqnarray}
{\boldsymbol  \theta}^{00}_{nm} &=& {\boldsymbol  \theta}^{00}-(n \pi\,sign(\theta^{00}_1)/2, m \pi\, sign(\theta^{00}_2)/2),\\
{\boldsymbol  \theta}^{ij}_{nm} &=& {\boldsymbol  \theta}^{00}_{nm}-(i \pi\,sign((\theta^{00}_{nm})_1), j \pi\, sign((\theta^{00}_{nm})_2)),
\end{eqnarray}
being $i,j,n,m \in\{0,1\}$.
Then, by denoting as $\widetilde{M}_h^{4h}({\boldsymbol \theta}^{00})$ the block-matrix representation of $M_h^{4h}$ on the Fourier space, the asymptotic three-grid convergence factor can be computed as the supreme of the spectral radii of the corresponding blocks
\begin{equation}
\rho_{3g}= \rho(M_h^{4h}) := \sup_{{\boldsymbol \theta}^{00}\in{\boldsymbol \Theta}_{4h}}\rho(\widetilde{M}_h^{4h}({\boldsymbol \theta}^{00})).
\end{equation}

In order to illustrate the suitability of the presented local Fourier analysis, next sections are devoted to the development of this analysis for two particular problems and to present some LFA results.

\section{Stokes problem}\label{sec:saddle}

As first model problem, we consider the Stokes equations. Essentially, the Stokes model represents slow viscous flow, it applies in particular to small-scale fluid movements in cells, lubrication, small particles, and porous-media-flow applications, in which we can ignore inertia forces.

\subsection{Description of the problem and the overlapping block smoother}\label{description_Stokes}

Stokes problem consists to find a velocity vector ${\mathbf u}: \Omega \rightarrow {\mathbb R}^2,$ ${\mathbf u} = (u,v),$ and a pressure field $p: \Omega \rightarrow {\mathbb R},$ such that

\begin{eqnarray}\label{Stokesp}
-\Delta {\mathbf u} + \nabla p &=& {\mathbf f}, \quad \hbox{in } \Omega, \nonumber \\
\nabla \cdot {\mathbf u} &=& 0, \quad \hbox{in } \Omega,\\
{\mathbf u} &=& {\mathbf 0}, \quad \hbox{on } \partial \Omega,\nonumber
\end{eqnarray}
where $\Omega$ is a bounded domain in ${\mathbb R}^2$.
To determine $p$ uniquely, making this problem well-posed, one needs to impose some additional condition, such as $$\int_{\Omega}p\, \hbox{d}{\mathbf x} = 0.$$\\
We are interested in using linear finite element methods to approximate the following variational formulation of problem (\ref{Stokesp}):\\
\noindent Find $(\mathbf{u}, p) \in
(H_0^1(\Omega))^2 \times L^2_0(\Omega)$ such that
\begin{eqnarray}
& & a({\mathbf u},{\mathbf v}) + (\nabla p,{\mathbf v}) =
({\mathbf f},{\mathbf v}),
 \quad \forall \, {\mathbf v} \in (H_0^1(\Omega))^2,  \\
& & (\nabla \cdot {\mathbf u},q) = 0,
\quad \forall q \in L^2_0(\Omega),
\end{eqnarray}
where $\displaystyle a({\mathbf u},{\mathbf v}) = \int_{\Omega} \nabla {\mathbf u} \cdot\nabla {\mathbf v}\, \hbox{d}{\mathbf x},$ and $L^2_0(\Omega)=\{ q \in L^2(\Omega)\;\vert \;\int_{\Omega} q \hbox{d} {\mathbf x} = 0\},$ that is,  $L^2(\Omega)$-functions which only differ by a constant are not distinguished.

Let ${\cal T}_h$ be a triangulation of $\Omega$ satisfying the usual admissibility assumption, and let ${\cal U}_h \subset (H_0^1(\Omega))^2,$ and ${\cal Q}_h \subset L^2_0(\Omega)$ be the corresponding spaces of piecewise linear functions on ${\cal T}_h.$ As the pair $({\cal U}_h, {\cal Q}_h)$ provides an unstable finite element scheme, we stabilize by modifying the discrete equations by introducing an additional term. To this end, the bilinear form on ${\cal Q}_h\times{\cal Q}_h,$ defined by
$$c(p_h,q_h) = \sum_{T \in {\cal T}_h} h_T^2 \int_{T} \nabla p_h \cdot \nabla q_h \hbox{d}{\mathbf x},$$
where $h_T$ denotes the diameter of the triangle $T,$ is considered. The stabilized discrete formulation of the Stokes problem in its weak form yields\\

\noindent Find $(\mathbf{u}_h, p_h) \in
{\cal U}_h \times {\cal Q}_h$ such that
\begin{eqnarray}
& & a({\mathbf u}_h,{\mathbf v}_h) + (\nabla p_h,{\mathbf v}_h) =
({\mathbf f},{\mathbf v}_h),
 \quad \forall \, {\mathbf v}_h \in {\cal U}_h,  \\
& & (\nabla \cdot {\mathbf u}_h,q_h) + \beta c(p_h, q_h) = 0,
\quad \forall q_h \in {\cal Q}_h,
\end{eqnarray}
where ${\mathbf u}_h = (u_h,v_h),$ the term $\beta c(p_h, q_h)$ refers to the stabilization of the problem, and $\beta > 0$ is a priori given parameter. The choice $\beta = 1/12$ appears to be optimal for linear elements~\cite{donea_huerta}, so it is used here as well.\\

\begin{figure}[htb]
\centering
\begin{tabular}{cc}
\begin{minipage}{5cm}
\includegraphics[scale=0.55]{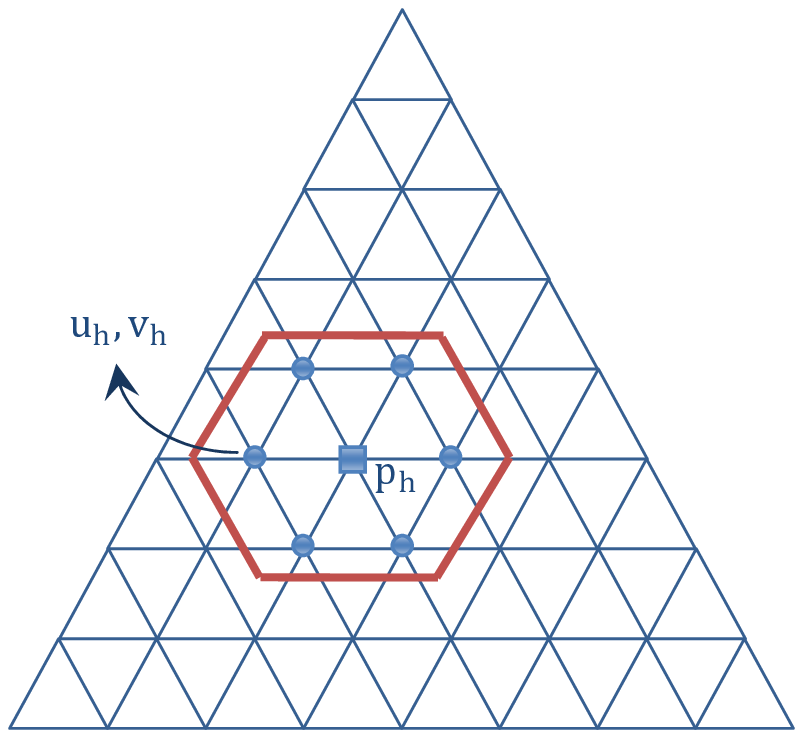}
\end{minipage}
&
\hspace{1cm}
\begin{minipage}{5cm}
\includegraphics[scale=0.55]{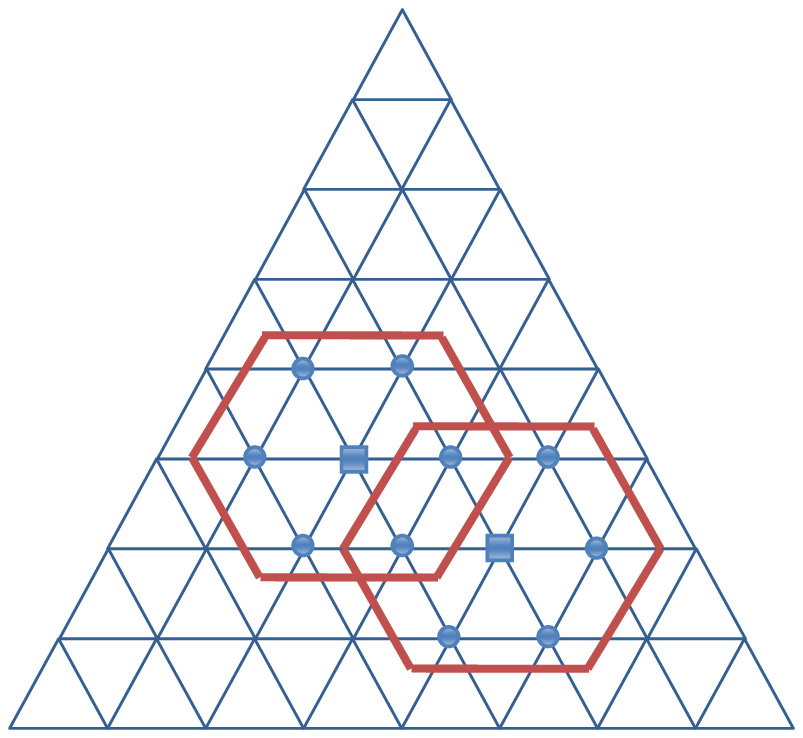}
\end{minipage}
\end{tabular}
\caption{Unknowns simultaneously updated in
point-wise box Gauss-Seidel smoother, and example of the overlapping between two arbitrary blocks. Circles denote velocity unknowns whereas the square refers to pressure degrees of freedom.}\label{vanka-smoothers}
\end{figure}

Overlapping block smoothers have been widely applied to saddle point problems, and in particular in the field of Computational Fluid Dynamics (CFD), see \cite{John, Tob, Turek1} for example. In this framework they are known as Vanka smoothers. Although the easiest overlapping block smoother corresponds to that updating all degrees of freedom corresponding to an element of the triangulation, it was shown that the use of this element-wise block smoother can be problematic for continuous pressure approximations, see~\cite{John_Matthies}. Then, as a good alternative the pressure-oriented Vanka smoother appears. Thus, more concretely, for Stokes equations, a suitable smoother is a point-wise box Gauss-Seidel iterative algorithm, which consists of simultaneously updating all unknowns appearing in the discrete divergence operator in the second equation of the system. This way of building the blocks is very common in box-relaxations used for Stokes and Navier-Stokes problems. The way in which the coupling between the pressure and velocity unknowns is treated, is very important for the convergence of the algorithm. The motivation which leads to this approach is that the incompressibility constraint ${\rm div} \, {\mathbf u},$ in general, is not satisfied locally in one element, whereas it can be fulfilled by increasing the number of involved velocities and considering the whole patch around one node.
This approach implies that $12$ unknowns corresponding to velocities and one pressure unknown (see left Figure \ref{vanka-smoothers}) are relaxed simultaneously and therefore, a $13\times13$ system has to be solved for each grid point. Notice that the discretization of divergence of ${\mathbf u}$ by linear finite element methods only involves the twelve velocity unknowns around the center node, and therefore to include the two velocity degrees of freedom located at this point is not necessary.
In this full variant (all the unknowns in the system are considered coupled), we iterate over all grid points in lexicographic order, and for each of them the corresponding box is solved.
For this version, the local system to solve for each box has the following form
\begin{equation}\label{point_vanka_system}
\left(
\begin{array}{cc}
{\mathbf M} & {\mathbf b} \\
{\mathbf b}^t & c
\end{array}
\right)
\left(
\begin{array}{c}
{\boldsymbol  \delta \mathbf u}\\
\delta p
\end{array}
\right)
=
\left(
\begin{array}{c}
{\mathbf r}^u\\
r^p
\end{array}
\right),
\end{equation}
being ${\mathbf M}$ a $12\times12$ matrix containing the coefficients associated with the velocity unknowns on the $12$ equations corresponding to these velocity unknowns, ${\mathbf b}$ the vector of the coefficients corresponding to the pressure unknown on the velocity equations, and $c$ the coefficient of the pressure point in its corresponding equation. In this way, ${\boldsymbol \delta \mathbf u}$ are the corrections of the velocity unknowns to be solved in the box, and $\delta p$ that of the pressure unknown. Finally, ${\mathbf r}^u$ and $r^p$ are the corresponding residuals of the unknowns appearing in the block.

The need of solving such systems makes these smoothers expensive. It is known that the smoothing is in general the most consuming part of the multigrid algorithm, and in this case it is even more noticeable due to the fact that the smoother consumes more than $70\%$ of the total computational time. A simplified variant can be considered by only coupling each velocity unknown with itself and the corresponding pressure unknown, in the previous system (\ref{point_vanka_system}). Thus, matrix ${\mathbf M}$ is transformed to a diagonal matrix ${\mathbf D},$ resulting the local systems in the following form
\begin{equation}\label{diag_vanka_system}
\left(
\begin{array}{cc}
{\mathbf D} & {\mathbf b} \\
{\mathbf b}^t & c
\end{array}
\right)
\left(
\begin{array}{c}
{\boldsymbol \delta \mathbf u}\\
\delta p
\end{array}
\right)
=
\left(
\begin{array}{c}
{\mathbf r}^u\\
r^p
\end{array}
\right).
\end{equation}
This box-type smoother is known as diagonal point-wise box-smoother, and the way in which the local systems are built allows an efficient implementation of such smoother. From system (\ref{diag_vanka_system}), the solutions can be computed as
\begin{equation*}
\delta p = \displaystyle \frac{{\mathbf b}^t{\mathbf D}^{-1}{\mathbf r}^u-r^p}{{\mathbf b}^t{\mathbf D}^{-1}{\mathbf b}-c},\qquad \quad {\boldsymbol \delta \mathbf u} = {\mathbf D}^{-1}({\mathbf r}^u- \delta p {\mathbf b}).
\end{equation*}
This resolution is based on performing some scalar products and componentwise products/additions between vectors, which results much cheaper than the solution of the full complete system (\ref{point_vanka_system}) by applying an LU decomposition or similar methods. Neglecting the time consumed in the construction of the LU decomposition, the diagonal version of the smoother requires roughly three times less operations than the full version, which makes the use of the former very appealing in practice.

\subsection{LFA for overlapping block smoother for Stokes problem}\label{LFA:stokes}

\begin{figure}[hbt]
\centering
\includegraphics[scale = 0.5]{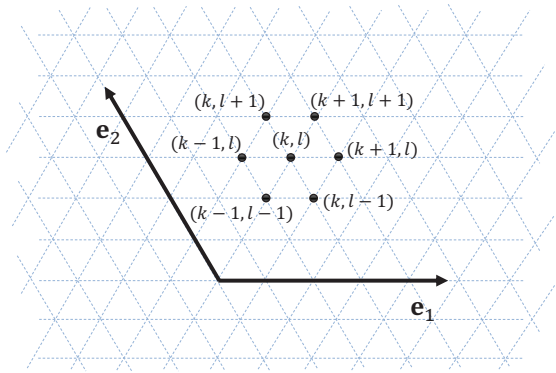}
\caption{Numbering of the nodes in an infinite regular triangular grid.}
\label{infinite_numbering}
\end{figure}

Next, we present some details of the local Fourier analysis for the point-wise box Gauss-Seidel in an infinite regular triangular grid.
With this purpose, we fix a suitable numbering of the grid-points. This is done by using a double index numeration, according to the unitary basis of ${\mathbb R}^2$ introduced in Section~\ref{sec:lfa}, see Figure~\ref{infinite_numbering}.

The considered point-wise box Gauss-Seidel simultaneously updates the unknown in node $(k,l)$ corresponding to the pressure $p_h,$ together with those unknowns corresponding to velocities $u_h$ and $v_h$ located in the six points around $(k,l),$ that is,  $(k+1,l),\, (k+1,l+1),\, (k,l+1),\, (k-1,l),\, (k-1,l-1),$ and $(k,l-1).$
Then, considering the equations corresponding to these unknowns, the resulting system to solve in terms of corrections and residuals reads
\begin{equation}\label{vanka_system_2}
A_h^B
\left(
\begin{array}{c}
\delta u_{k+1,l}\\
\delta u_{k+1,l+1}\\
\delta u_{k,l+1}\\
\delta u_{k-1,l}\\
\delta u_{k-1,l-1}\\
\delta u_{k,l-1}\\
\delta v_{k+1,l}\\
\delta v_{k+1,l+1}\\
\delta v_{k,l+1}\\
\delta v_{k-1,l}\\
\delta v_{k-1,l-1}\\
\delta v_{k,l-1}\\
\delta p_{k,l}
\end{array}
\right)=
\left(
\begin{array}{c}
r^u_{k+1,l}\\
r^u_{k+1,l+1}\\
r^u_{k,l+1}\\
r^u_{k-1,l}\\
r^u_{k-1,l-1}\\
r^u_{k,l-1}\\
r^v_{k+1,l}\\
r^v_{k+1,l+1}\\
r^v_{k,l+1}\\
r^v_{k-1,l}\\
r^v_{k-1,l-1}\\
r^v_{k,l-1}\\
r^p_{k,l}
\end{array}
\right).
\end{equation}
In the case of the diagonal version, matrix $A_h^B$ is much more sparse than for the full-Vanka smoother, since each velocity unknown is only coupled with the pressure variable appearing in the system.
Notice that this system corresponds to system~\eqref{vanka_system} for the case $m=3$ and with \\
\begin{center}
\begin{tabular}{ccc}
${\boldsymbol \delta}{\mathbf u}^1 = \left(\begin{array}{c}\delta u_{k+1,l}\\ \delta u_{k+1,l+1}\\ \delta u_{k,l+1}\\ \delta u_{k-1,l}\\ \delta u_{k-1,l-1}\\ \delta u_{k,l-1}\end{array}\right),$
&
${\boldsymbol \delta}{\mathbf u}^2 = \left(\begin{array}{c}\delta v_{k+1,l}\\ \delta v_{k+1,l+1}\\ \delta v_{k,l+1}\\ \delta v_{k-1,l}\\ \delta v_{k-1,l-1}\\ \delta v_{k,l-1}\end{array}\right),$
&
${\boldsymbol \delta}{\mathbf u}^3 =  \delta p_{k,l}.$
\end{tabular}
\end{center}
Taking into account that velocity unknowns are updated up to $s=6$ times and pressure unknowns are only updated once during each complete relaxation step, it is fulfilled that
\begin{eqnarray*}
\delta u_{k+kk,l+ll} \!\!\!&=&\!\!\! u_h^{j+(n+1)/6}({\mathbf x}_{k+kk,l+ll})-u_h^{j+n/6}({\mathbf x}_{k+kk,l+ll}),\\ [1ex]
\delta v_{k+kk,l+ll} \!\!\!&=&\!\!\! v_h^{j+(n+1)/6}({\mathbf x}_{k+kk,l+ll})-v_h^{j+n/6}({\mathbf x}_{k+kk,l+ll}),\\ [1ex]
\delta p_{k,l} \!\!\!&=&\!\!\! p_h^{j+1}({\mathbf x}_{k,l})-p_h^{j}({\mathbf x}_{k,l}),
\end{eqnarray*}
with $kk,ll = -1,0,1$, excluding the case $kk = ll = 0$, and $n=0,\ldots,5.$
As commented in Section~\ref{sec:lfa}, to find the relation between the initial and the final errors, we must write system (\ref{vanka_system_2}) in terms of the errors, taking into account the intermediate errors. Let us denote as $e^j_{u,h}({\mathbf x}),$ $e^j_{v,h}({\mathbf x}),$ $e^j_{p,h}({\mathbf x}),$ the initial error for each variable $u_h,$ $v_h,$ and $p_h,$ at j-iteration; $e^{j+1}_{u,h}({\mathbf x}),$ $e^{j+1}_{v,h}({\mathbf x}),$ $e^{j+1}_{p,h}({\mathbf x}),$ the corresponding initial errors at $(j+1)-$iteration; and $e^{j+n/6}_{u,h}({\mathbf x}),$ $e^{j+n/6}_{v,h}({\mathbf x}),\; n=1,\ldots,5,$ the intermediate errors for the velocities, obtained after the unknown has been updated $n-$times in the current iteration. Next, following the assumption done in~(\ref{errors_fourier}), we consider
\begin{equation}\label{errors}
\begin{array}{lll}
e_{u,h}^j({\mathbf x}) = \alpha_u^{(0)}({\boldsymbol \theta})\, e^{\imath\boldsymbol \theta {\mathbf x}}, & e_{v,h}^j({\mathbf x}) = \alpha_v^{(0)}({\boldsymbol \theta})\, e^{\imath\boldsymbol \theta {\mathbf x}}, & e_{p,h}^j({\mathbf x}) = \alpha_p^{(0)}({\boldsymbol \theta})\, e^{\imath\boldsymbol \theta {\mathbf x}},\\ [1ex]
e_{u,h}^{j+1}({\mathbf x}) = \alpha_u^{(6)}({\boldsymbol \theta})\, e^{\imath\boldsymbol \theta {\mathbf x}}, & e_{v,h}^{j+1}({\mathbf x}) = \alpha_v^{(6)}({\boldsymbol \theta})\, e^{\imath\boldsymbol \theta {\mathbf x}}, & e_{p,h}^{j+1}({\mathbf x}) = \alpha_p^{(6)}({\boldsymbol \theta})\, e^{\imath\boldsymbol \theta {\mathbf x}},\\ [1ex]
e_{u,h}^{j+n/6}({\mathbf x}) = \alpha_u^{(n)}({\boldsymbol \theta})\, e^{\imath\boldsymbol \theta {\mathbf x}}, & e_{v,h}^{j+n/6}({\mathbf x}) = \alpha_v^{(n)}({\boldsymbol \theta})\, e^{\imath\boldsymbol \theta {\mathbf x}}, &
n = 1,\ldots,5.
\end{array}
\end{equation}
Notice that the Fourier coefficient of the fully corrected error for the pressure has been denoted as $\alpha_p^{(6)},$ instead of $\alpha_p^{(1)},$ although the pressure unknown is fully corrected when it has been updated only once. This notation is chosen only to match with that of the velocities.\\
Taking into account that the corrections can be written in terms of the errors, i.e.
\begin{eqnarray*}
\delta u_{k+kk,l+ll} \!\!\!&=&\!\!\! e_{u,h}^{j+(n+1)/6}({\mathbf x}_{k+kk,l+ll})-e_{u,h}^{j+n/6}({\mathbf x}_{k+kk,l+ll}),\\ [1ex]
\delta v_{k+kk,l+ll} \!\!\!&=&\!\!\! e_{v,h}^{j+(n+1)/6}({\mathbf x}_{k+kk,l+ll})-e_{v,h}^{j+n/6}({\mathbf x}_{k+kk,l+ll}),\\ [1ex]
\delta p_{k,l} \!\!\!&=&\!\!\! e_{p,h}^{j+1}({\mathbf x}_{k,l})-e_{p,h}^{j}({\mathbf x}_{k,l}),
\end{eqnarray*}
which are given by expressions in~(\ref{errors}), and that before updating the block of unknowns corresponding to node $(k,l),$ the variables involved in equations in~(\ref{vanka_system}) are in the states showing in Table~\ref{updated_variables},
the corrections appearing in system~(\ref{vanka_system_2}) can be written as
\vspace{-0.3cm}
\begin{equation*}
\left(
\begin{array}{c}
\delta u_{k+1,l}\\
\delta u_{k+1,l+1}\\
\delta u_{k,l+1}\\
\delta u_{k-1,l}\\
\delta u_{k-1,l-1}\\
\delta u_{k,l-1}\\
\delta v_{k+1,l}\\
\delta v_{k+1,l+1}\\
\delta v_{k,l+1}\\
\delta v_{k-1,l}\\
\delta v_{k-1,l-1}\\
\delta v_{k,l-1}\\
\delta p_{k,l}
\end{array}
\right)=
\left(
\begin{array}{c}
(\alpha_u^{(3)}(\boldsymbol \theta)-\alpha_u^{(2)}(\boldsymbol \theta))e^{\imath\theta_1}\\
(\alpha_u^{(1)}(\boldsymbol \theta)-\alpha_u^{(0)}(\boldsymbol \theta))e^{\imath\theta_1}\,e^{\imath\theta_2}\\
(\alpha_u^{(2)}(\boldsymbol \theta)-\alpha_u^{(1)}(\boldsymbol \theta))e^{\imath\theta_2}\\
(\alpha_u^{(4)}(\boldsymbol \theta)-\alpha_u^{(3)}(\boldsymbol \theta))e^{-\imath\theta_1}\\
(\alpha_u^{(6)}(\boldsymbol \theta)-\alpha_u^{(5)}(\boldsymbol \theta))e^{-\imath\theta_1}\,e^{-\imath\theta_2}\\
(\alpha_u^{(5)}(\boldsymbol \theta)-\alpha_u^{(4)}(\boldsymbol \theta))e^{-\imath\theta_2}\\
(\alpha_v^{(3)}(\boldsymbol \theta)-\alpha_v^{(2)}(\boldsymbol \theta))e^{\imath\theta_1}\\
(\alpha_v^{(1)}(\boldsymbol \theta)-\alpha_v^{(0)}(\boldsymbol \theta))e^{\imath\theta_1}\,e^{\imath\theta_2}\\
(\alpha_v^{(2)}(\boldsymbol \theta)-\alpha_v^{(1)}(\boldsymbol \theta))e^{\imath\theta_2}\\
(\alpha_v^{(4)}(\boldsymbol \theta)-\alpha_v^{(3)}(\boldsymbol \theta))e^{-\imath\theta_1}\\
(\alpha_v^{(6)}(\boldsymbol \theta)-\alpha_v^{(5)}(\boldsymbol \theta))e^{-\imath\theta_1}\,e^{-\imath\theta_2}\\
(\alpha_v^{(5)}(\boldsymbol \theta)-\alpha_v^{(4)}(\boldsymbol \theta))e^{-\imath\theta_2}\\
(\alpha_p^{(6)}(\boldsymbol \theta)-\alpha_p^{(0)}(\boldsymbol \theta))\\
\end{array}
\right) e^{\imath{\boldsymbol \theta}{\mathbf x}}.\\
\end{equation*}
\begin{table}[htb]
\begin{center}
\begin{tabular}{|c|c|c|c|c|c|c|}
\hline
0 & 1 & 2 & 3 & 4 & 5 & 6 \\
\hline
\fbox{$p_{k,l}$}     & \fbox{$u_{k,l+1}$} & $u_{k-1,l+1}$ & \fbox{$u_{k-1,l}$} & \fbox{$u_{k,l-1}$}   & \fbox{$u_{k-1,l-1}$} & $u_{k,l-2}$   \\
\fbox{$u_{k+1,l+1}$} & \fbox{$v_{k,l+1}$} & $v_{k-1,l+1}$ & \fbox{$v_{k-1,l}$} & \fbox{$v_{k,l-1}$}   & \fbox{$v_{k-1,l-1}$} & $v_{k,l-2}$   \\
\fbox{$v_{k+1,l+1}$} &             & \fbox{$u_{k+1,l}$}   &             & $u_{k+1,l-1}$ &               & $u_{k-1,l-2}$ \\
$u_{k,l+2}$   &             & \fbox{$v_{k+1,l}$}   &             & $v_{k+1,l-1}$ &               & $v_{k-1,l-2}$ \\
$v_{k,l+2}$   &             & $u_{k+2,l}$   &             & $u_{k-2,l}$   &               & $u_{k-2,l-2}$ \\
$u_{k+1,l+2}$ &             & $v_{k+2,l}$   &             & $v_{k-2,l}$   &               & $v_{k-2,l-2}$ \\
$v_{k+1,l+2}$ &             &               &             &               &               & $u_{k-2,l-1}$ \\
$u_{k+2,l+1}$ &             &               &             &               &               & $v_{k-2,l-1}$ \\
$v_{k+2,l+1}$ &             &               &             &               &               &               \\
$u_{k+2,l+2}$ &             &               &             &               &               &               \\
$v_{k+2,l+2}$ &             &               &             &               &               &               \\
\hline
\end{tabular}
\end{center}
\caption{Number of times that each variable has been relaxed before updating the block of unknowns corresponding to node $(k,l).$ The variables involved in such a block are framed to highlight those that are also updated in the iteration.}
\label{updated_variables}
\end{table}
In the same way, the residuals composing the right-hand side in (\ref{vanka_system_2}) can also be expressed in terms of coefficients $\alpha_u^{(n)}(\boldsymbol\theta),$ $\alpha_v^{(n)}(\boldsymbol\theta),$ and $\alpha_p^{(n)}(\boldsymbol\theta).$ Thus, the resulting system can be rewritten as a system of equations for the updated Fourier coefficients: $\alpha_u^{(n)}(\boldsymbol\theta),$ $\alpha_v^{(n)}(\boldsymbol\theta),\; n=1,\ldots,6,$ and $\alpha_p^{(6)}(\boldsymbol\theta),$ leaving in the right-hand side the terms corresponding to the non-updated components, $\alpha_u^{(0)}(\boldsymbol\theta),$ $\alpha_v^{(0)}(\boldsymbol\theta),$ and $\alpha_p^{(0)}(\boldsymbol\theta),$ as explained in~\eqref{systemPQ}. In this particular case, $P$ is a $(13\times13)-$matrix, $Q$ is a $(13\times3)-$matrix, and then matrix $P^{-1}\,Q$, appearing in~\eqref{systemPQ_2},  is $(13\times3).$ Finally, from the obtained system, a relation between the initial and fully corrected errors can be obtained by considering the final $(3\times3)-$block of matrix $(P^{-1}\,Q),$ i.e. $(P^{-1}\,Q)_{\{11:13,1:3\}},$ which consists of the last three rows, resulting in the following
\begin{equation*}
\left( \begin{array}{c}\alpha_u^{(6)}(\boldsymbol \theta) \\ \alpha_v^{(6)}(\boldsymbol \theta) \\ \alpha_p^{(6)}(\boldsymbol \theta) \end{array} \right) = \;(P^{-1}\,Q)_{\{11:13,1:3\}}\; \left( \begin{array}{c} \alpha_u^{(0)}(\boldsymbol \theta) \\ \alpha_v^{(0)}(\boldsymbol \theta) \\ \alpha_p^{(0)}(\boldsymbol \theta) \end{array}\right) = \widetilde{{\mathbf S}_h}(\boldsymbol \theta) \left( \begin{array}{c} \alpha_u^{(0)}(\boldsymbol \theta) \\ \alpha_v^{(0)}(\boldsymbol \theta) \\ \alpha_p^{(0)}(\boldsymbol \theta) \end{array}\right).
\end{equation*}
From $\widetilde{{\mathbf S}_h}(\boldsymbol \theta)$ the smoothing as well as the k-grid Fourier analysis for this smoother can then be performed as standard local Fourier analysis on triangular grids (see Section~\ref{sec:lfa}).

\subsection{Numerical results}\label{sec:Stokes_results}

This section is devoted to show some results from the presented local Fourier analysis in order to demonstrate its accuracy and utility. With this purpose, we compare the two- and three-grid convergence factors provided by LFA to real asymptotic convergence factors. These latter are computed by using a grid obtained after performing nine refinement levels to an equilateral triangle with unit edge length. Zero right-hand side and a random initial guess are considered to avoid round-off errors.\\

In Table~\ref{tabla_2g_Stokes}, we show smoothing $\mu$ and two-grid convergence factors $\rho_{2g}$ obtained by LFA, together with the corresponding asymptotic convergence factors, $\rho_h$, computed by using the multigrid version of the algorithm. Results for different numbers of smoothing steps $\nu$ are displayed for both full- and diagonal-versions of the smoother. First, we observe that the analysis predicts very accurately the convergence factors computationally obtained by using a $W-$cycle. Also we can see that the convergence provided by the diagonal version is very similar to that of the full smoother. This makes preferable the use of this latter, due to its lower computational cost.
\begin{table}[htb]
\begin{center}
\begin{tabular}{|c|c|c|c|c|c|c|}
\cline{2-7}
\multicolumn{1}{c}{} & \multicolumn{3}{|c|}{Full-Vanka} & \multicolumn{3}{|c|}{Diagonal-Vanka}  \\
\hline
$\nu= \nu_1+\nu_2$ & \multirow{2}{*}{$\mu^{\nu}$} & \multirow{2}{*}{$\rho_{2g}$} & \multirow{2}{*}{$\rho_h$} & \multirow{2}{*}{$\mu^{\nu}$} & \multirow{2}{*}{$\rho_{2g}$} & \multirow{2}{*}{$\rho_h$} \\
& & & & & & \\
\hline
$1$ & $0.51$ & $0.64$ & $0.63$ & $0.55$ & $0.69$ & $0.69$ \\
\hline
$2$ & $0.26$ & $0.34$ & $0.34$ & $0.31$ & $0.31$ & $0.30$ \\
\hline
$3$ & $0.14$ & $0.18$ & $0.17$ & $0.17$ & $0.19$ & $0.19$ \\
\hline
$4$ & $0.07$ & $0.13$ & $0.12$ & $0.09$ & $0.15$ & $0.15$ \\
\hline
$5$ & $0.04$ & $0.10$ & $0.10$ & $0.05$ & $0.13$ & $0.12$ \\
\hline
\end{tabular}
\end{center}
\caption{Smoothing, $\mu^{\nu}$, and two-grid convergence factors, $\rho_{2g}$, predicted by LFA together with real asymptotic convergence factors, $\rho_{h}$, obtained by using a W-cycle, for different numbers of smoothing steps $\nu = \nu_1+\nu_2$ and both considered overlapping block smoothers.}
\label{tabla_2g_Stokes}
\end{table}
In many situations, the use of a relaxation parameter in the multigrid smoothing can accelerate the convergence dramatically. However, there is no rule in order to choose this parameter. Here, we use local Fourier analysis in order to choose the optimal relaxation parameters, that is, the relaxation parameters that provide the best two-grid convergence factor. With this purpose, we perform a loop over the possible values of $\omega_u$ and $\omega_p$, that is, the relaxation parameters for velocities and for pressure, respectively. For example, if three smoothing steps are considered, $\omega_u=1.05$ and $\omega_p= 0.6$ provide an optimal two-grid convergence factor of $\rho_{2g}=0.087$ for the multigrid based on the diagonal overlapping block smoother. This factor predicted by LFA matches perfectly with the asymptotic convergence factor computationally obtained by using a multigrid version of a $W-$cycle. However, as well-known this type of multigrid cycle is expensive and we are interested in the use of $V-$cycles, if possible. To see if this is suitable, we have developed a three-grid local Fourier analysis for the diagonal version of the smoother. This three-grid analysis allows to see the different behavior of $V-$ and $W-$cycles. Then, in Table~\ref{tabla_3g_Stokes}, we show the three-grid convergence factors predicted by the LFA for a $V-$cycle by using different numbers of pre- and post-smoothing steps. Together, we present the computationally obtained asymptotic convergence factors by using a multigrid V-cycle and a fine grid with nine refinement levels. In the case of only one pre-smoothing step, the multigrid shows divergence (although when a three-grid method the factor matches that predicted by the LFA, of course). However, in the rest of the tests we can observe that the computationally obtained results match perfectly those predicted by the analysis.
\begin{table}[htb]
\begin{center}
\begin{tabular}{|c|c|c|}
\hline
\multirow{2}{*}{$(\nu_1,\nu_2)$} & \multirow{2}{*}{$\rho_{3g}$} & \multirow{2}{*}{$\rho_h$} \\
& & \\
\hline
$(1,0)$ & $0.68$ & div \\
\hline
$(1,1)$ & $0.31$ & $0.31$ \\
\hline
$(2,1)$ & $0.28$ & $0.28$ \\
\hline
$(2,2)$ & $0.24$ & $0.23$ \\
\hline
\end{tabular}
\end{center}
\caption{Three-grid convergence factors, $\rho_{3g}$, predicted by LFA together with real asymptotic convergence factors, $\rho_{h}$, for different numbers of pre- and post-smoothing steps $(\nu_1,\nu_2)$ by using a $V-$cycle and the diagonal smoother.}
\label{tabla_3g_Stokes}
\end{table}
Moreover, as previously done with the two-grid LFA, we can find the relaxation parameters providing the optimal three-grid convergence factor for a V-cycle, by using the proposed three-grid analysis. With this tool, we find that an optimal three-grid convergence factor of $\rho_{3g}=0.097$ can be obtained by using $\omega_u = 0.95$ and $\omega_p =0.6$ and the diagonal overlapping block smoother.
In order to see if this strategy is useful from the practical point of view, we are going to consider a well-known benchmark problem for testing new
computational algorithms. We consider a cavity flow problem consisting of a two-dimensional triangle enclosure with the upper side translating with uniform velocity. The boundary conditions are no slip on sides of the triangle moving with a velocity of constant magnitude, and on fixed sides the velocity is zero. This problem is solved on a unit equilateral triangle.
To perform this simulation, we have applied the multigrid algorithm based on the diagonal version of the smoother by using $V-$ and $W-$cycles for different numbers of refinement levels. For both cycles, we have used the corresponding relaxation parameters $\omega_u$ and $\omega_p$ to obtain optimal three-grid convergence.
These results are shown in Table~\ref{tabla_levels_Stokes}, where the number of iterations necessary to reduce the initial residual in a factor of $10^{-10}$ are displayed. In the case of $W-$cycles, we observe that the convergence is independent on the number of levels used in the multigrid method. When $V-$cycles are considered, we can observe a slight increase in the number of iterations needed for the desired convergence. Despite this, the use of $V-$cycles seems to be suitable from the practical point of view, since the computational cost is much lower and the convergence is still very satisfactory. \\
\begin{table}[htb]
\begin{center}
\begin{tabular}{|c|c|c|c|c|c|c|c|}
\cline{2-8}
\multicolumn{1}{c|}{} &\multirow{2}{*}{$4$ lev.} & \multirow{2}{*}{$5$ lev.} & \multirow{2}{*}{$6$ lev.} & \multirow{2}{*}{$7$ lev.} & \multirow{2}{*}{$8$ lev.} & \multirow{2}{*}{$9$ lev.} & \multirow{2}{*}{$10$ lev.} \\
\multicolumn{1}{c|}{} & & & & & & & \\
\hline
$W-$cycle & $10$ & $10$ & $10$ & $10$ & $10$ & $10$ & $10$  \\
\hline
$V-$cycle & $10$ & $10$ & $11$ & $11$ & $12$ & $13$ & $14$  \\
\hline
\end{tabular}
\end{center}
\caption{Number of iterations necessary to reduce the initial residual in a factor of $10^{-10},$ for different numbers of refinement levels, by using $W-$ and $V-$cycles and the diagonal version of the smoother.}
\label{tabla_levels_Stokes}
\end{table}

\section{Vector model problem}\label{sec:vector}

In this section, we use a simple vector model problem in order to describe in detail the local Fourier analysis for overlapping block smoothers for edge-based discretizations of vector problems. Mainly as a result of its good computational properties, edge-based discretizations have emerged widely in the simulation of many real applications, including electromagnetic field computations.

\subsection{Description of the problem and the overlapping block smoother}\label{description_vector}

In particular, we consider the following model problem:
\begin{eqnarray}\label{rot_rot_u1}
\hbox{\textbf{curl}}\, \hbox{rot}\, {\mathbf u} + \kappa\,{\mathbf u} = {\mathbf f}, \;\; \hbox{in } \; \Omega,\label{rot_rot_u}\\
{\mathbf u}\times {\mathbf n} = 0, \; \; \hbox{on} \; \partial \Omega,\label{rot_rot_u1}
\end{eqnarray}
where $\Omega$ is an open domain in ${\mathbb R}^2,$ ${\mathbf f}\in (L^2(\Omega))^2,$
and ${\mathbf n}$ denotes the outward unit normal vector along the boundary $\partial \Omega.$ We assume vanishing tangential components on the boundary~(\ref{rot_rot_u1}).
As usual, if we define the Hilbert space~(see~\cite{Girault_Raviart})
\begin{equation}\label{H0_curl}
{\mathbf H_0}(\hbox{rot},\Omega) := \{{\mathbf u}\in (L^2(\Omega))^2 \, | \, \hbox{rot}\,{\mathbf u}\in L^2(\Omega), \, {\mathbf u }\times{\mathbf n} = 0\; \hbox{on} \; \partial\Omega\},
\end{equation}
the weak formulation of problem~\eqref{rot_rot_u} reads\\
\noindent Find $\mathbf{u}\in {\mathbf H_0}(\hbox{rot},\Omega)$ such that
\begin{equation}\label{weak_form}
a({\mathbf u},{\mathbf v}) = ({\mathbf f},{\mathbf v}), \quad \forall \,{\mathbf v}\in {\mathbf H_0}(\hbox{rot},\Omega)
\end{equation}
where $\displaystyle a({\mathbf u},{\mathbf v}) = \int_{\Omega}(\hbox{rot}\, {\mathbf u})\,(\hbox{rot}\, {\mathbf v})\,\hbox{d{\textbf x}} + \kappa\int_{\Omega}{\mathbf u}\cdot {\mathbf v}\,\hbox{d{\textbf x}}$

To discretize this variational problem, first we consider a triangulation of $\Omega$ satisfying the usual admissibility assumption, ${\mathcal T}_h$. Canonical finite elements for the approximation of ${\mathbf H_0}(\hbox{rot},\Omega)$ are the one by N\'{e}d\'{e}lec~\cite{Nedelec}. Here, low-order N\'{e}d\'{e}lec's edge elements are considered. That is, we consider the family of vectors which on every element $K$ of the triangulation ${\cal T}_h$ are linear in each component and tangential continuous across the element edges, that is
\begin{equation}\label{Vh}
{\mathbf V}_h = \{{\mathbf v}_h\in {\mathbf H}_0(\hbox{rot},\Omega) \, | \, {\mathbf v}_h|_K = \left[\begin{array}{c}a_1 \\ a_2 \end{array}\right]+b \left[\begin{array}{c}y \\ -x \end{array}\right],\, \forall K\in {\cal T}_h \}.
\end{equation}
With this choice, the finite element approximation of~\eqref{weak_form} reads\\
\noindent Find $\mathbf{u}_h\in {\mathbf V}_h$ such that
\begin{equation}\label{discrete_weak_form}
a({\mathbf u}_h,{\mathbf v}_h) = ({\mathbf f}_h,{\mathbf v}_h), \quad \forall \,{\mathbf v_h}\in {\mathbf V}_h.
\end{equation}

\begin{figure}[htb]
\begin{center}
\begin{tabular}{cc}
\begin{minipage}{5cm}
\includegraphics[scale=0.5]{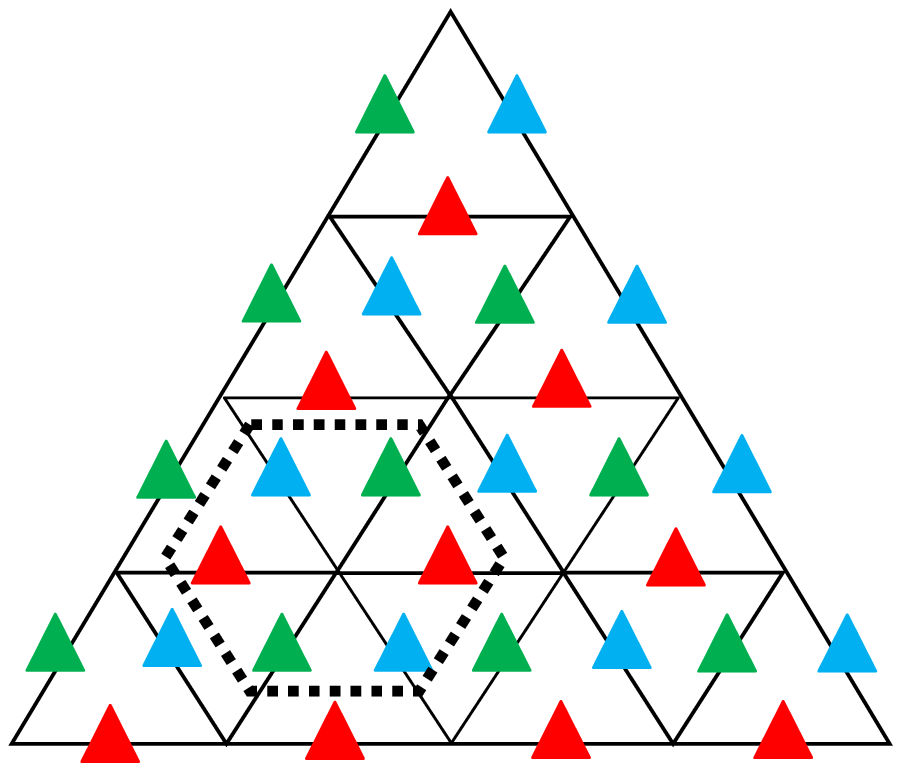}
\end{minipage}
&
\hspace{1cm}
\begin{minipage}{5cm}
\includegraphics[scale=0.5]{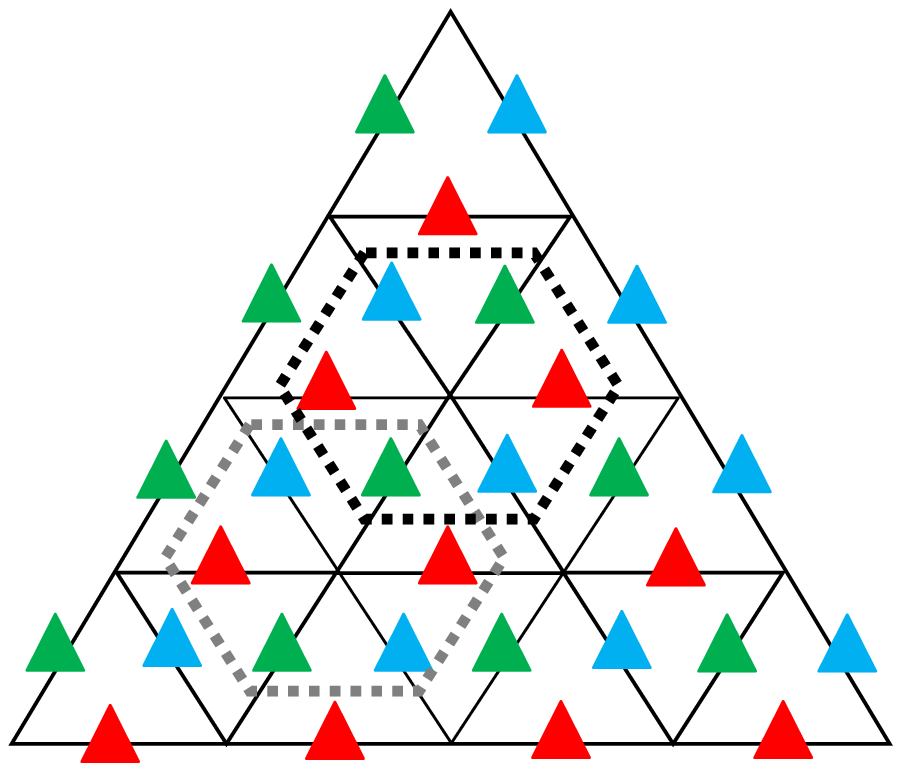}
\end{minipage}
\\
(a) & (b)
\end{tabular}
\end{center}
\caption{(a) Unknowns simultaneously updated in
point-wise overlapping block smoother for the vector problem, and (b) overlapping of the blocks.}\label{point-vanka-smoother}
\end{figure}

Some efforts to design multigrid solution schemes for edge-based discretizations have been carried out recently, see for example~\cite{Hiptmair, HiptmairXu, Hu_et_al, Kolev1, Kolev2, Reitzinger}. It is important the choice of a suitable smoother, and following the suggestion in some of these works (e.g.~\cite{Arnold}), we are going to consider an appropriate Schwarz--type smoother. Within this relaxation procedure, for every vertex of the triangular grid, all edges connected with such a vertex are smoothed simultaneously. Therefore, a $(6\times 6)-$system has to be solved for each vertex in the grid. This set of unknowns is shown in Figure~\ref{point-vanka-smoother} (a), and in the same figure (Figure~\ref{point-vanka-smoother} (b)) we can observe the overlapping between these blocks of unknowns that gives rise to the Schwarz character of the smoother.

\subsection{Local Fourier analysis for the overlapping block smoother}

In order to perform a local Fourier analysis for this type of edge-based discretizations, we have to take into account some of their special characteristics.
In this type of discretizations, there are unknowns located at different types of grid-points, and therefore the stencils defining the discrete operator on each point-type involve different surrounding unknowns, what makes that the discrete operator is not defined in the same way at every grid-point. This fact has to be considered in the performance of the analysis for this type of discretizations.

\begin{figure}[htb]
\begin{center}
\begin{tabular}{cc}
\includegraphics[scale = 0.7]{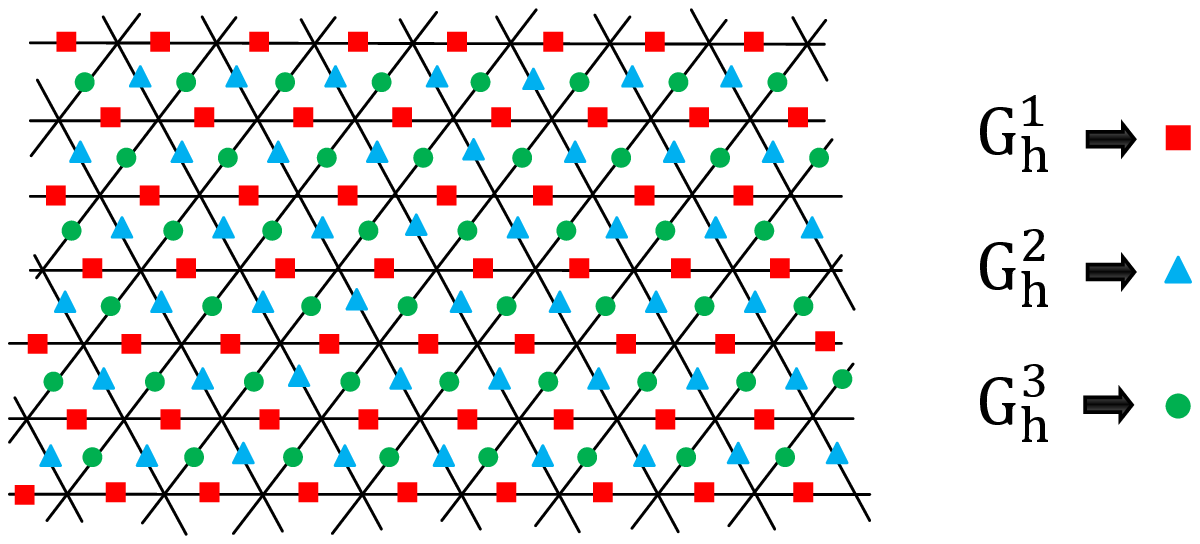}
& \hspace{0.3cm}
\includegraphics[scale = 0.5]{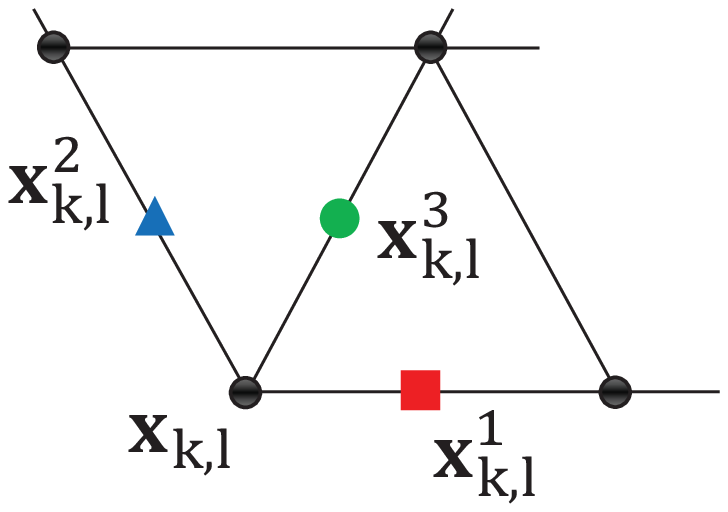}
\\
(a) & \hspace{0.3cm} (b)
\end{tabular}
\end{center}
\caption{(a) Infinite grid, composed of three different infinite subgrids. (b) Location of the different unknowns, and correspondence to the three different subgrids, together with the local numbering for the grid-points.}
\label{infinite_grid2}
\end{figure}

\begin{figure}[htb]
\begin{center}
\begin{tabular}{ccc}
\includegraphics[scale = 0.55]{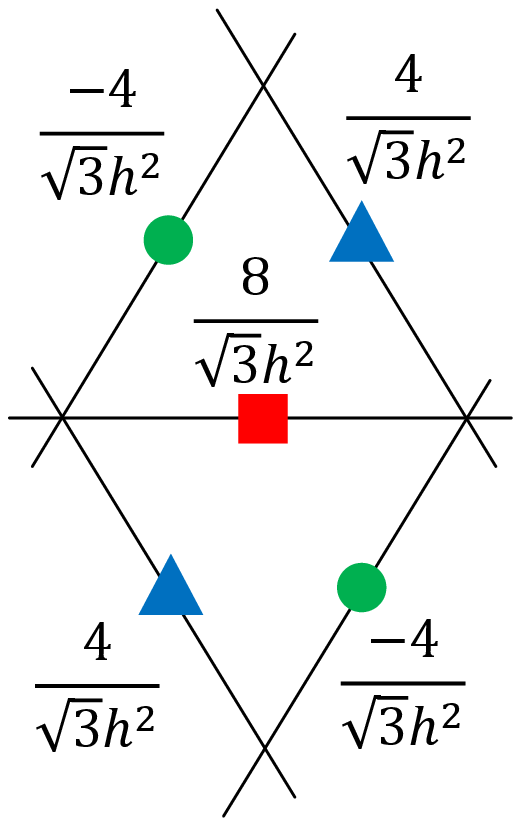}
\hspace{0.5cm}
&
\includegraphics[scale = 0.55]{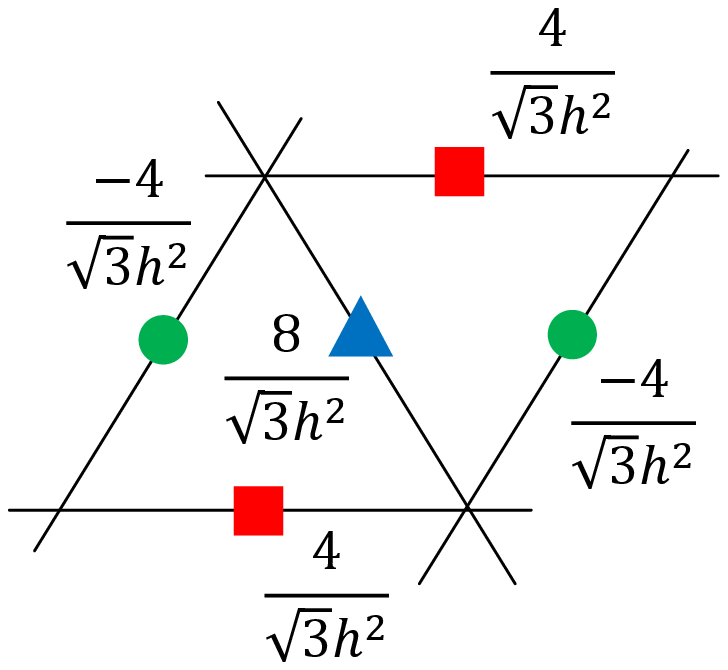}
\hspace{0.5cm}
&
\includegraphics[scale = 0.55]{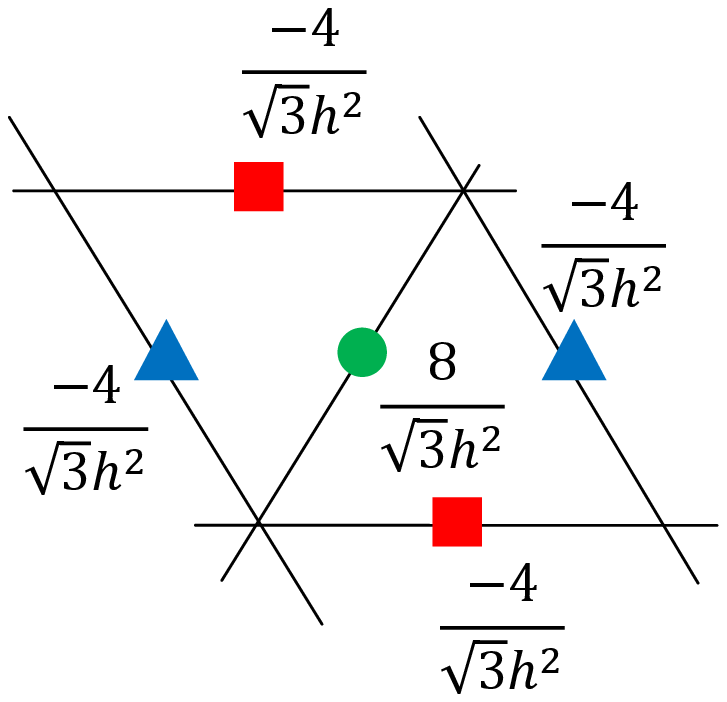}
\\
(a)\hspace{0.5cm} &  (b)\hspace{0.5cm} & (c)
\end{tabular}
\end{center}
\caption{Stencils for operator ${\bf curl} \, {\rm rot}$ obtained by the lowest order N\'ed\'elec FEM for equilateral triangular grids, corresponding to grid-points located at (a) $G_h^1$, (b) $G_h^2$, and (c) $G_h^3$.}
\label{Nedelec_stencils}
\end{figure}

First of all, as an infinite grid has to be considered, the original grid is extended to the following infinite grid $G_h = \displaystyle\bigcup_{i=1}^{3} G_h^i,$ see Figure~\ref{infinite_grid2} (a), given as the union of three different subgrids
\begin{equation}\label{infinite_grid}
G_h^i:=\{{\mathbf x}_{k,l}^i = ((k+\delta_1^i)\,h_1\,{\mathbf e}_1, (l+\delta_2^i)\,h_2\,{\mathbf e}_2)\,\vert\, k,l\in{\mathbb Z}\},
\end{equation}
where $$(\delta_1^i,\delta_2^i)=\left\{\begin{array}{ll}(1/2,0),&\hbox{if}\, i=1,\\ (0,1/2),&\hbox{if}\, i=2,\\(1/2,1/2),&\hbox{if}\, i=3.\end{array}\right.$$
Notice that each of these subgrids $G_h^j$ are associated with a type of edge, that is to each one of the three different orientations, as we can see in Figure~\ref{infinite_grid2} (b).
In order to extend the definition of the discrete operator to the infinite grid $G_h$, we have to take into account that the resulting equations at grid-points on $G_h^1$, $G_h^2$ and $G_h^3$ are different. Thus, we can define the application of the discrete operator to a grid function $u_h$ on $G_h$ in the following way:
\begin{eqnarray}
\hspace{-0.5cm}A_h\,u_h ({\mathbf x}) &=& \displaystyle\left\{\begin{array}{l}\displaystyle A_h^{11} u_h({\mathbf x_{k,l}^1}) + A_h^{12} u_h({\mathbf x_{k,l}^2})+ A_h^{13} u_h({\mathbf x_{k,l}^3}), \;{\mathbf x} = {\mathbf x}_{k,l}^1 \in G_h^1 \\ [0.5ex]
A_h^{21} u_h({\mathbf x_{k,l}^1}) + A_h^{22} u_h({\mathbf x_{k,l}^2})+ A_h^{23} u_h({\mathbf x_{k,l}^3}), \;{\mathbf x} = {\mathbf x}_{k,l}^2 \in G_h^2\\ [0.5ex]
A_h^{31} u_h({\mathbf x_{k,l}^1}) + A_h^{32} u_h({\mathbf x_{k,l}^2})+ A_h^{33} u_h({\mathbf x_{k,l}^3}), \;{\mathbf x} = {\mathbf x}_{k,l}^3 \in G_h^3
\end{array}\right.
\label{discrete_operator}
\\ [1ex]
&=&\displaystyle\left\{\begin{array}{l}
\displaystyle\sum_{r = 1}^{3}\Big(\sum_{(kk,ll)\in I^{1 r}} s_{kk,ll}^{1 r} u_h({\mathbf x_{k+kk,l+ll}^{r}})\Big), \;{\mathbf x} = {\mathbf x}_{k,l}^1 \in G_h^1, \\ [0.8ex]
\displaystyle\sum_{r = 1}^{3}\Big(\sum_{(kk,ll)\in I^{2 r}} s_{kk,ll}^{2 r} u_h({\mathbf x_{k+kk,l+ll}^{r}})\Big), \;{\mathbf x} = {\mathbf x}_{k,l}^2 \in G_h^2,\\ [0.8ex]
\displaystyle\sum_{r = 1}^{3}\Big(\sum_{(kk,ll)\in I^{3 r}} s_{kk,ll}^{3 r} u_h({\mathbf x_{k+kk,l+ll}^{r}})\Big), \;{\mathbf x} = {\mathbf x}_{k,l}^3 \in G_h^3,
\end{array}\right.\nonumber
\end{eqnarray}
where values $s_{kk,ll}^{i r}$ are the coefficients corresponding to the stencil of discrete operator $A_h^{i r}$ which gives the relation that exists in the equation between one unknown in $G_h^i$ and the unknowns in $G_h^{r}$. Subsets $I^{i r}$ give the connections of a grid-point located at $G_h^i$ with those in $G_h^{r}$.
In order to illustrate the definition of the discrete operator, we are going to consider the lowest order N\'ed\'elec finite element discretization of operator ${\bf curl} \, {\rm rot}$ on equilateral triangular grids.
Writing $A_h = N_h + M_h$, where $N_h$ and $M_h$ represent the matrix corresponding to the ${\bf curl} \, {\rm rot}$ operator and the mass matrix, respectively, the stencils corresponding to $N_h$ for the three different grid-points are given in Figure~\ref{Nedelec_stencils}.
Besides, we can write these stencils with the notation previously given in the following way:\\

\begin{tabular}{ccc}
\hspace{-0.5cm}$N_h^{11} = \frac{8}{\sqrt{3}h^2}\left[\begin{array}{ccc}0 & 0 & 0 \\ 0 & 1 & 0 \\ 0 & 0 & 0\end{array}\right]$,
&
$N_h^{12} = \frac{4}{\sqrt{3}h^2}\left[\begin{array}{ccc}0 & 0 & 0 \\ 0 & 0 & 1 \\ 0 & 1 & 0\end{array}\right]$,
&
$N_h^{13} = \frac{-4}{\sqrt{3}h^2}\left[\begin{array}{ccc}0 & 0 & 0 \\ 0 & 1 & 0 \\ 0 & 1 & 0\end{array}\right]$,
\\
\\
\hspace{-0.5cm}$N_h^{21} = \frac{4}{\sqrt{3}h^2}\left[\begin{array}{ccc}0 & 1 & 0 \\ 1 & 0 & 0 \\ 0 & 0 & 0\end{array}\right]$,
&
$N_h^{22} = \frac{8}{\sqrt{3}h^2}\left[\begin{array}{ccc}0 & 0 & 0 \\ 0 & 1 & 0 \\ 0 & 0 & 0\end{array}\right]$,
&
$N_h^{23} = \frac{-4}{\sqrt{3}h^2}\left[\begin{array}{ccc}0 & 0 & 0 \\ 1 & 1 & 0 \\ 0 & 0 & 0\end{array}\right]$,
\\
\\
\hspace{-0.5cm}$N_h^{31} = \frac{-4}{\sqrt{3}h^2}\left[\begin{array}{ccc}0 & 1 & 0 \\ 0 & 1 & 0 \\ 0 & 0 & 0\end{array}\right]$,
&
$N_h^{32} = \frac{-4}{\sqrt{3}h^2}\left[\begin{array}{ccc}0 & 0 & 0 \\ 0 & 1 & 1 \\ 0 & 0 & 0\end{array}\right]$,
&
$N_h^{33} = \frac{8}{\sqrt{3}h^2}\left[\begin{array}{ccc}0 & 0 & 0 \\ 0 & 1 & 0 \\ 0 & 0 & 0\end{array}\right]$.
\\
\\
\end{tabular}

Following the description of the LFA for overlapping block smoothers in Section~\ref{sec:lfa}, next we are going to explain the corresponding analysis for our concrete problem. With this aim, we denote $u_h^i = u_h|_{G_h^i},$ for $i=1,2,3$.
In this case, system~\eqref{vanka_system} is fulfilled with $m =3$ and \\
\begin{center}
\begin{tabular}{ccc}
${\boldsymbol \delta}{\mathbf u}^1 = \left(\begin{array}{c}\delta u^1_{k,l}\\ \delta u^1_{k-1,l}\end{array}\right),$
&
${\boldsymbol \delta}{\mathbf u}^2 = \left(\begin{array}{c}\delta u^2_{k,l}\\ \delta u^2_{k,l-1}\end{array}\right),$
&
${\boldsymbol \delta}{\mathbf u}^3 = \left(\begin{array}{c}\delta u^3_{k,l}\\ \delta u^3_{k-1,l-1}\end{array}\right).$
\end{tabular}
\end{center}

\noindent Notice that each unknown is updated up to $s=2$ times per relaxation step, and it is fulfilled that
$$\delta u^i_{k,l} = u_h^{i,j+(n+1)/2}({\mathbf x}_{k,l}^i)-u_h^{i,j+n/2}({\mathbf x}_{k,l}^i),\; i=1,2,3,\; n =0,1.$$ Again, denoting the initial, final and intermediate errors as $(e_h^i)^j$, $(e_h^i)^{j+1}$ and $(e_h^i)^{j+1/2}$, respectively, it holds
$$
\begin{array}{lcr}
e_h^{i,j} = \alpha_i^{(0)}({\boldsymbol \theta})e^{\imath{\boldsymbol \theta}{\mathbf x}},\quad & e_h^{i,j+1} = \alpha_i^{(2)}({\boldsymbol \theta})e^{\imath{\boldsymbol \theta}{\mathbf x}},\quad & e_h^{i,j+1/2}=\alpha_i^{(1)}({\boldsymbol \theta})e^{\imath{\boldsymbol \theta}{\mathbf x}}\quad,
\end{array}
$$
and the corrections can be written in terms of the errors in the following way
$$\delta u^i_{k,l} = e_h^{i,j+(n+1)/2}({\mathbf x}_{k,l}^i)-e_h^{i,j+n/2}({\mathbf x}_{k,l}^i).$$
\begin{table}[htb]
\renewcommand
{\arraystretch}{1.5}
\begin{center}
\begin{tabular}{|c|c|c|}
\hline
0 & 1 & 2 \\
\hline
\fbox{$u^1_{k,l}$}     & \fbox{$u^1_{k-1,l}$}   & $u^1_{k-1,l-1}$ \\
$u^1_{k,l+1}$   & \fbox{$u^2_{k,l-1}$}   & $u^2_{k-1,l-1}$ \\
\fbox{$u^2_{k,l}$}     & $u^3_{k-1,l}$   &                 \\
$u^2_{k+1,l}$   & \fbox{$u^3_{k-1,l-1}$} &                 \\
\fbox{$u^3_{k,l}$}     & $u^3_{k,l-1}$   &                 \\ [0.7ex]
\hline
\end{tabular}
\end{center}
\caption{Number of times that each variable has been relaxed before updating the block of unknowns corresponding to vertex $(k,l).$ The variables involved in such a block are framed to highlight those that are also updated in the iteration.}
\label{updated_variables_vector}
\end{table}
Taking into account that before updating the block of unknowns around vertex $(k,l)$, the variables involved in~\eqref{vanka_system} are in the states shown in Table~\ref{updated_variables_vector}, the corrections in~\eqref{vanka_system} can be written as
\vspace{-0.3cm}
\begin{equation*}
\left(
\begin{array}{c}
\delta u^1_{k,l}\\
\delta u^1_{k-1,l}\\
\delta u^2_{k,l}\\
\delta u^2_{k,l-1}\\
\delta u^3_{k,l}\\
\delta u^3_{k-1,l-1}
\end{array}
\right)=
\left(
\begin{array}{c}
(\alpha_1^{(1)}(\boldsymbol \theta)-\alpha_1^{(0)}(\boldsymbol \theta))\\
(\alpha_1^{(2)}(\boldsymbol \theta)-\alpha_1^{(1)}(\boldsymbol \theta))e^{-\imath\theta_1}\\
(\alpha_2^{(1)}(\boldsymbol \theta)-\alpha_2^{(0)}(\boldsymbol \theta))\\
(\alpha_2^{(2)}(\boldsymbol \theta)-\alpha_2^{(1)}(\boldsymbol \theta))e^{-\imath\theta_2}\\
(\alpha_3^{(1)}(\boldsymbol \theta)-\alpha_3^{(0)}(\boldsymbol \theta))\\
(\alpha_3^{(2)}(\boldsymbol \theta)-\alpha_3^{(1)}(\boldsymbol \theta))e^{-\imath\theta_1}\, e^{-\imath\theta_2} \\
\end{array}
\right) e^{\imath{\boldsymbol \theta}{\mathbf x}}.\\
\end{equation*}

Similarly, the residuals appearing in the right-hand side of the system can also be expressed in terms of coefficients $\alpha_i^{(n)}({\boldsymbol \theta})$, and therefore we can write a system for the updated Fourier coefficients: $\alpha_i^{(n)}({\boldsymbol \theta})$, $n = 1,2,$ that is,\\
\begin{tabular}{ccc}
\begin{minipage}{6cm}
\begin{equation*}
P\;\left( \begin{array}{c}\alpha_1^{(1)}(\boldsymbol \theta)\\ \alpha_2^{(1)}(\boldsymbol \theta)\\  \alpha_3^{(1)}(\boldsymbol \theta)\\ \alpha_1^{(2)}(\boldsymbol \theta)\\ \alpha_2^{(2)}(\boldsymbol \theta)\\  \alpha_3^{(2)}(\boldsymbol \theta)\\ \end{array} \right) = \;Q\; \left( \begin{array}{c} \alpha_1^{(0)}(\boldsymbol \theta) \\ \alpha_2^{(0)}(\boldsymbol \theta) \\ \alpha_3^{(0)}(\boldsymbol \theta) \end{array}\right)
\end{equation*}
\end{minipage}
&
\begin{minipage}{1cm}
\begin{center}
\includegraphics[width=0.8cm, height=1.2cm]{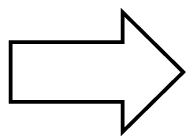}
\end{center}
\end{minipage}
&
\begin{minipage}{6cm}
\begin{equation*}
\left( \begin{array}{c}\alpha_1^{(1)}(\boldsymbol \theta)\\ \alpha_2^{(1)}(\boldsymbol \theta)\\  \alpha_3^{(1)}(\boldsymbol \theta)\\ \alpha_1^{(2)}(\boldsymbol \theta)\\ \alpha_2^{(2)}(\boldsymbol \theta)\\  \alpha_3^{(2)}(\boldsymbol \theta)\\ \end{array} \right) = \;(P^{-1}\,Q)\; \left( \begin{array}{c} \alpha_1^{(0)}(\boldsymbol \theta) \\ \alpha_2^{(0)}(\boldsymbol \theta) \\ \alpha_3^{(0)}(\boldsymbol \theta) \end{array}\right),
\end{equation*}
\end{minipage}\\
\vspace{0.2cm}
\end{tabular}\\
where $P$ is a $(6\times 6)-$matrix, $Q$ is a $(6\times 3)-$matrix, and $P^{-1}Q$ is a $(6\times 3)-$matrix. By considering the final $(3\times 3)-$block of matrix $P^{-1}Q$, we can obtain the relation between the initial and fully corrected errors as follows
\begin{equation*}
\left( \begin{array}{c}\alpha_1^{(2)}(\boldsymbol \theta) \\ \alpha_2^{(2)}(\boldsymbol \theta) \\ \alpha_3^{(2)}(\boldsymbol \theta) \end{array} \right) = \;(P^{-1}\,Q)_{\{4:6,1:3\}}\; \left( \begin{array}{c} \alpha_1^{(0)}(\boldsymbol \theta) \\ \alpha_2^{(0)}(\boldsymbol \theta) \\ \alpha_3^{(0)}(\boldsymbol \theta) \end{array}\right) = \widetilde{{\mathbf S}_h}(\boldsymbol \theta) \left( \begin{array}{c} \alpha_1^{(0)}(\boldsymbol \theta) \\ \alpha_2^{(0)}(\boldsymbol \theta) \\ \alpha_3^{(0)}(\boldsymbol \theta) \end{array}\right).
\end{equation*}
In this way, we obtain the Fourier representation of the overlapping block smoother, $\widetilde{{\mathbf S}_h}(\boldsymbol \theta)$, and we can perform the smoothing as well as the k-grid Fourier analysis.

\subsection{Numerical results}\label{sec:vector_results}

In this section we present some results from the developed local Fourier analysis compared to real asymptotic convergence factors. These latter are obtained by using a grid obtained after performing nine refinement levels to an equilateral triangle with unit edge length. Zero right-hand side and a random initial guess are considered to avoid round-off errors.
First, we consider problem~\eqref{rot_rot_u}-\eqref{rot_rot_u1} with $\kappa = 1$. In Table~\ref{tabla_vector}, we show smoothing $\mu$ and three-grid convergence factors $\rho_{3g}$ obtained by LFA, together with the corresponding asymptotic convergence factors. Results for a number of combinations of pre- and post-smoothing steps as well as for different type of cycles (V- and W-cycles) are displayed.
\begin{table}[htb]
\begin{center}
\begin{tabular}{|c|c|c|c|c|c|c|}
\cline{4-7}
\multicolumn{3}{c}{} & \multicolumn{2}{|c|}{V-cycle} & \multicolumn{2}{|c|}{W-cycle}  \\
\hline
\multirow{2}{*}{$\nu$} & \multirow{2}{*}{$(\nu_1,\nu_2)$} & \multirow{2}{*}{$\mu^{\nu}$} & \multirow{2}{*}{$\rho_{3g}$} & \multirow{2}{*}{$\rho_h$} & \multirow{2}{*}{$\rho_{3g}$} & \multirow{2}{*}{$\rho_h$} \\
& & & & & & \\
\hline
\multirow{2}{*}{$1$} & $(1,0)$ & \multirow{2}{*}{$0.46$} & $0.34$ & $0.33$ & \multirow{2}{*}{$0.33$} & \multirow{2}{*}{$0.33$} \\
& $(0,1)$ & & $0.34$ & $0.33$ & & \\
\hline
\multirow{3}{*}{$2$} & $(1,1)$ & \multirow{3}{*}{$0.21$} & $0.13$ & $0.13$ & \multirow{3}{*}{$0.12$} & \multirow{3}{*}{$0.12$} \\
& $(2,0)$ & & $0.13$ & $0.13$ & & \\
& $(0,2)$ & & $0.13$ & $0.13$ & & \\
\hline
\multirow{4}{*}{$3$} & $(2,1)$ & \multirow{4}{*}{$0.09$} & $0.07$ & $0.07$ & \multirow{4}{*}{$0.07$} & \multirow{4}{*}{$0.07$} \\
& $(1,2)$ & & $0.07$ & $0.07$ & & \\
& $(3,0)$ & & $0.07$ & $0.07$ & & \\
& $(0,3)$ & & $0.07$ & $0.07$ & & \\
\hline
\end{tabular}
\end{center}
\caption{Smoothing, $\mu^{\nu}$, and three-grid convergence factors, $\rho_{3g}$, predicted by LFA together with real asymptotic convergence factors, $\rho_{h}$, for different numbers of smoothing steps $\nu = \nu_1+\nu_2$.}
\label{tabla_vector}
\end{table}
First of all, it is observed that the LFA results predict with very high accuracy the values of the real convergence factors. In this case, from the analysis it is observed that the behavior of V- and W-cycles is very similar and that no difference is seen when different combinations of pre- and post-smoothing steps are used. Therefore, we can deduce that only two smoothing steps are necessary to obtain a convergence factor around $0.1$, and that V-cycles are preferred against W-cycles since they provide the same convergence factors but with a lower computational cost, keeping the parallel features of the multigrid algorithm.
We have observed from the analysis that the use of a relaxation parameter doesn't give rise to an improvement in the behavior of the method.

Finally, we show the robust behavior of the multigrid method regarding coefficient $\kappa$ and the number of unknowns. A $V(1,1)-$cycle has been chosen due to the satisfactory factors predicted by the LFA when $\kappa = 1$. In Table~\ref{table_different_kappa}, we present for different values of $\kappa$ and various refinement levels, the number of iterations necessary to reduce the initial residual in a factor of $10^{-10}$.
\begin{table}[htb]
\begin{center}
\begin{tabular}{|l|c|c|c|c|}
\hline
& 6 levels & 7 levels & 8 levels & 9 levels \\
\hline
$\kappa = 1$       & 11 & 11 & 11 & 11 \\
$\kappa = 10^{-2}$ & 11 & 11 & 11 & 11 \\
$\kappa = 10^{-4}$ & 11 & 11 & 11 & 11 \\
$\kappa = 10^{-8}$ & 11 & 11 & 11 & 11 \\
\hline
\end{tabular}
\end{center}
\caption{$V(1,1)-$cycle convergence for different values of $\kappa$ and several numbers of refinement levels.}
\label{table_different_kappa}
\end{table}
We observe that the number of iterations is independent with respect to the value of $\kappa$ and the number of refinement levels.

\section{Conclusions}\label{sec:conclusions}

In this work, a general framework to perform a local Fourier analysis for overlapping block smoothers on triangular grids is presented. A detailed description of the computation of the Fourier representation of the smoothing operator  has been done for the general case, and the way to perform a two- and three-grid local Fourier analysis has been also explained. In order to clarify the development of this analysis, two particular examples have been considered. First, we have presented the local Fourier analysis of the so-called Vanka smoother for the multigrid solution of the Stokes equations. For this particular problem, we have seen that the performance of the diagonal version of the smoother is comparable to that of the full Vanka relaxation and we have seen its suitability by displaying results from the two- and three-grid LFA, as well as, some results with real multigrid cycles. The performance of V-cycles is also compared to that of the W-cycles, yielding to more efficient algorithms. On the other hand, a vector model problem based on the
curl-curl operator is also considered. For an edge-based discretization of this problem, the analysis of a suitable overlapping block smoother has been performed. Also, results from the LFA and from real multigrid cycles are presented to demonstrate the usefulness of the analysis. Summarizing, one can say that the general framework for the LFA for overlapping block smoothers can be used for different relaxations of this kind and for different types of discretizations.

\section*{Acknowledgements}

This research has been partially supported by FEDER/MCYT Projects MTM2013-40842-P and the DGA (Grupo consolidado PDIE).






\bibliographystyle{model1b-num-names}
\bibliography{bib_Vanka}

\end{document}